\magnification 1200 
 \pageno=0
 \hsize=140mm
 \vsize=200mm
 \hoffset=-4mm
 \voffset=-1mm
 \pretolerance=500
 \tolerance=1000
 \brokenpenalty=5000
 \frenchspacing

\font\sc=cmcsc10
\let\eps =\varepsilon
\def\reel{{\rm I}\!{\rm R}}
\def\rond{{\scriptstyle\circ}}
\def\N#1{\left\Vert#1\right\Vert}
\def\abs#1{\left\vert#1\right\vert}
\let\ph=\varphi
\def\frac#1#2{{#1\over#2}}
\def\comp{\;{}^{{}_\vert}\!\!\!{\rm C}}
\def\nat{{{\rm I}\!{\rm N}}}
\def\cqfd{\unskip\kern 6pt\penalty 500
\raise -2pt\hbox{\vrule\vbox to 10pt{\hrule width 4pt\vfill
\hrule}\vrule}\par}
\def\tvi{\vrule height 12pt depth 5pt width 0pt}
\def\Ssi{\quad\Longleftrightarrow\quad}
\def\bg{\bigskip\goodbreak}
\def\adh#1{\overline{\strut#1}}
\def\ref#1&#2&#3&#4&#5\par{\par{\leftskip = 5em{\noindent
\kern-5em\vbox{\hrule height0pt depth0pt width
5em\hbox{\bf[\kern2pt#1\unskip\kern2pt]\enspace}}\kern0pt}
{\sc\ignorespaces#2\unskip},\
{\rm\ignorespaces#3\unskip}\
{\sl\ignorespaces#4\unskip\/}\
{\rm\ignorespaces#5\unskip}\par}}
\def\sn{\nobreak\smallskip}

 \font\eightrm=cmr8
\font\eightbf=cmbx8
\font\eighti=cmmi8
\font\eightrm=cmr8
\font\eightronde=cmmi10 at 14pt
\font\eightsc=cmcsc8
\font\eightsl=cmsl8
\font\eightsy=cmsy8

\def\eightpoint{\textfont0=\eightrm\def\rm{\fam0\eightrm}
\textfont1=\eighti\scriptfont1=\sixi \textfont2=\eightsy
\def\bf{\eightbf}
\def\sc{\eightsc}
\def\sl{\eightsl}
\normalbaselineskip=9pt \normalbaselines
\eightrm}

\font\ronde=cmmi10 at 16pt
\font\sixi=cmmi6
\font\sronde=cmmi10 at 12pt

\def\A{{\cal A}}
\def\B{{\cal B}}
\def\C{{\cal C}}
\def\E{{\cal E}}
\def\erA{{\hbox{\eightronde a}}}
\def\K{{\cal K}}
\def\M{{\cal M}}
\def\P{{\cal P}}
\def\rA{{\hbox{\ronde a}}}
\def\srA{{\hbox{\sronde a}}}
\def\U{{\cal U}}
\def\V{{\cal V}}
\def\1{{\bf 1}}

\def\ds{\displaystyle}
\def\equival{\mathrel{\Longleftrightarrow}}
\def\implies{\mathrel{\Longrightarrow}}
\def\proof{\noindent{\bf Proof: }}
\def\Sim{\mathop{\,\sim\,}}
\def\span{{\rm span\,}}
\def\t{\tilde}
\def\To{\mathop{\longrightarrow}}
\def\Tr{{\rm Tr}}

\centerline{\bf On Subspaces of Non-commutative $L_p$-Spaces}
\vskip 0.5cm
\centerline{\bf Yves Raynaud, Quanhua Xu}
\vskip 2cm

{\eightpoint {\bf Abstract:}  We study some structural aspects of
the subspaces of the non-commutative (Haagerup) $ L_p$-spaces
associated with a general (non necessarily semi-finite) von
Neumann algebra $\erA$. If a subspace $ X$ of $ L_p(\erA)$
contains uniformly the spaces $\ell_p^n$, $ n\geq 1$, it contains
an almost isometric, almost 1-complemented copy of $\ell_p$. If $
X$ contains uniformly the finite dimensional Schatten classes $
S_p^n$, it contains their $\ell_p$-direct sum too. We obtain a
version of the classical Kadec-Pe\l czy\'nski dichotomy theorem
for $ L_p$-spaces, $ p\geq 2$. We also give operator space
versions of these results. The proofs are based on previous
structural results on the ultrapowers of $ L_p(\erA)$, together
with a careful analysis of the elements of an ultrapower $
L_p(\erA)_\U$ which are disjoint from the subspace $ L_p(\erA)$.
These techniques permit to recover a recent result of N.
Randrianantoanina concerning a Subsequence Splitting Lemma for the
general non-commutative $ L_p$ spaces. Various notions of $
p$-equiintegrability are studied (one of which is equivalent to
Randrianantoanina's one) and some results obtained by  Haagerup,
Rosenthal and Sukochev for $ L_p$-spaces based on finite von
Neumann algebras concerning subspaces of $ L_p(\erA)$ containing $
\ell_p$ are extended to the general case.}

 {\eightpoint \parindent= 0pt\footnote{} {2000 {\sl Mathematics
Subject Classification:} { Primary: 46B20, 46L52.
Secondary: 47M07\hfill\break
 {\sl Key words:} Non-commutative $L_p$-spaces, Ultrapowers, Schatten classes,
Equiintegrability}}}

\vskip 2cm

{\bf Contents:}

\bigskip

 \item{0.} Introduction
 \item{1.} Preliminaries
 \item{2.} Elements in $L_p(\rA)_\U$ which are
disjoint from $L_p(\rA)$
 \item{3.} Embedding of $\ell_p$-sums of finite
dimensional spaces
 \item{4.} Equiintegrability and the Subsequence
Splitting Lemma
 \item{5.} Subspaces containing $\ell_p$
 \item{6.} Operator space version
 \item{}   Appendix:  Equality case in non-commutative
Clarkson inequality

\vfill\eject

\centerline{\bf 0. Introduction}

\medskip

 Since several years the study of non-commutative $L_p$-spaces
has incited new interest because of their close relations with the
new and rapidly developing Operator Space Theory and
Non-commutative Probability Theory. It is known now that
non-commutative integration is a fundamental tool in both latter
theories. Conversely, results and problems from these theories
permit to gain new insight into the theory of non-commutative
$L_p$-spaces and at the same time pose new problems in the frame
of this theory: see for instance the recent works [EJR], [Ju1-3],
[JX], [NO], [O], [PX], [R1-4].

The starting point of the present work is a problem arising from
the theory of ${\cal OL}_p$-spaces, which was initiated by Effros
and Ruan ([ER1]) and developed in the recent paper [JNRX] (see
also [JOR], [NO]). To explain this, we first recall the famous
Kadec-Pe\l czy\'nski dichotomy theorem, which states that every
closed subspace of $L_p(0,1)$, $2<p<\infty$, either is isomorphic
to a Hilbert space or contains a subspace which is isomorphic to
$\ell_p$ and complemented in $L_p(0,1)$. This theorem plays an
important role in the classical theory of ${\cal L}_p$-spaces. The
class of ${\cal OL}_p$-spaces is an analog for the category of
operator spaces of the class of ${\cal L}_p$-spaces in the
category of Banach spaces; when going back to the Banach space
category by forgetting the matricial structure, the class of
${\cal OL}_p$-spaces gives rise to a still new class of Banach
spaces, which could be called ``non-commutative ${\cal
L}_p$-spaces'': $X$ belongs to this class if for some $\lambda$,
every finite dimensional subspace of $X$ is contained in another
subspace, which is $\lambda$-isomorphic to a finite dimensional
non-commutative $L_p$-space. It is then natural to look for a
non-commutative version of the Kadec-Pe\l czy\'nski dichotomy. A
version  of the usual Kadec-Pe\l czy\'nski's dichotomy in a
non-commutative setting exists already in the literature, with
exactly the same  statement; it was proved indeed in [S] for
non-commutative $L_p$-spaces based on  finite von Neumann algebras
(under an equivalent form),
 and in [R2] for semi-finite ones.

This version however does not help at all in developing the theory
of ${\cal OL}_p$-spaces (the conclusion it gives is ``too
commutative'' in a certain sense). A stronger, still very
hypothetical statement would be preferable in this direction:

\smallskip\noindent{\it A closed subspace of a non-commutative
$L_p$-space ($2<p<\infty$) should either be embeddable into a
\underbar{commutative} $L_p$-space or contain a copy of the
$p$-direct sum $K_p=(\mathop{\bigoplus}\limits_{n\geq 1}S_p^n)_p$
of the finite dimensional $p$-Schatten classes.}\smallskip

 A step towards this direction is made in the present paper, namely the
following theorem, which is one of our main results:

\proclaim Theorem {0.1}. Let $\rA$ be a von Neumann algebra (non
necessarily semi-finite), $0<p<\infty$, $p\neq 2$ and $X$ a closed
subspace of $L_p(\rA)$. Assume that for some constant $\lambda
\geq 1$, and for every $n\geq 1$, $X$ contains a subspace
$\lambda$-isomorphic to the space $S_p^n$ (resp. and
$\mu$-complemented in $L_p(\rA)$ -- in this case we suppose $p\geq
1$).  Then for every $\eps>0$, $X$ contains a subspace
$(\lambda+\eps)$-isomorphic to $K_p$ (resp. and
$(\lambda\mu+\eps)$-complemented in $L_p(\rA)$).

This result has a forerunner in the case $\rA=B(\ell_2)$ (then
$L_p(\rA)$ is the usual Schatten class $S_p$), which was obtained
by Arazy and Lindenstrauss in [ArL]. Their proof, which relies on
a careful analysis of the local structure of $S_p$ together with a
clever use of Ramsey's theorem, can be extended to some special
cases of Theorem 0.1 (e.g. when $\rA$ is finite and $p>2$) but we
hardly imagine how to adapt it to the general situation described
in Theorem~0.1. Our proof of Theorem 0.1 heavily depends on
ultrapower techniques, using the fact, proved in [Ra],  that the
class of non-commutative $L_p$-spaces is closed under ultrapowers.
In fact, we will see that the subspace of $X$ isomorphic to $K_p$
obtained in Theorem 0.1 is built over a sequence of subspaces
isomorphic to the $S_p^n$'s and ``almost disjoint''. This approach
to Theorem 0.1 also allows us to extend to all von Neumann
algebras the first non-commutative version of  the Kadec-Pe\l
czy\'nski dichotomy mentioned previously, which remained an open
question in the non semi-finite case. More precisely, we have:

\proclaim Theorem {0.2}. Let $\rA$ be a von Neumann algebra,
$2<p<\infty$ and $X$ a closed subspace of $L_p(\rA)$. Then either
$X$ is isomorphic to a Hilbert space and complemented in
$L_p(\rA)$ or $X$ contains a subspace isomorphic to $\ell_p$ and
complemented in $L_p(\rA)$.

This paper is organized as follows. In section 1 we recall some
necessary preliminaries on non-commutative $L_p$-spaces and their
ultrapowers. The non-commutative $L_p$-spaces we consider are
those constructed by Haagerup [H]. Contrary to the class of
``usual'' $L_p$-spaces associated with a  normal faithful
semi-finite trace, the class of Haagerup $L_p$-spaces is closed
under ultraproducts ([Ra]). The main tools of the paper are
developed in section 2, where we show how to push disjoint
elements in an ultrapower of $L_p(\rA)$ down to disjoint elements
of $L_p(\rA)$. Theorem 0.1 above will be proved in Section 3. In
fact, we shall prove a more general result by replacing the spaces
$S_p^n$ by a sequence of finite dimensional spaces. Section 4 is
devoted to the equiintegrability in $L_p(\rA)$. We give a rather
complete study of the $p$-equiintegrable subsets of $L_p(\rA)$.
Our techniques permit us to easily recover the Subsequence
Splitting Lemma proved by N. Randrianantoanina [R3]. In section $5$
we characterize the subspaces of $L_p(\rA)$ which contain a
subspace isomorphic to $\ell_p$. As a corollary, we  get Theorem
0.2. Such characterizations are classical in the commutative case,
and were recently proved for spaces associated with finite or
semifinite von Neumann algebras in [HRS], [R1] and [SX]. The last
section aims at extending the previous results to the operator
space setting. There we get the operator space versions of
Theorems 0.1 and 0.2. We also add an appendix whose result
determines when equality occurs in the non-commutative Clarkson
inequality. This result improves a previous theorem due to H.
Kosaki [Ko2] and implies a characterization of isometric 2-dimensional
$\ell_p$-subspaces of $L_p(\rA)$ which is repeatedly used in the paper.

\smallskip

The main results of this paper were announced in the Note [RaX].

\bigskip\centerline{\bf 1. Preliminaries}

\medskip
 This section contains notations, most notions and basic facts
 necessary to the whole paper. For clarity we divide it into three
 subsections.
\medskip

\noindent{\bf Non-commutative $L_p$-spaces}

\medskip There are several equivalent constructions of non-commutative
$L_p$-spaces associated with a von Neumann algebra (c.f., e.g.
[AM], [H], [Hi], [I], [Ko1], [Te2]). We shall use in this paper
Haagerup's construction, which we recall briefly now (see [Te1]
for a precise introduction to the subject). Let $\rA$ be a von
Neumann algebra. For $0<p<\infty$, the spaces $L_p(\rA)$ are
constructed as spaces of measurable operators relative not to
$\rA$ but to a certain semi-finite super von Neumann algebra  of
$\rA$, namely, the crossed product of $\rA$ by one of its modular
automorphism groups. Let $\M$ be the crossed product of $\rA$ by
the modular automorphism group $(\sigma_t)_{t\in\reel}$ of a fixed
normal faithful semifinite weight $w$ on $\rA$ (see [KaR], II.13).
Let $(\theta_s)$ be the dual automorphism group on $\M$. It is
well known that $\rA$ is a von Neumann subalgebra of $\M$ and that
the position of $\rA$ in $\M$ is determined by the group
$(\theta_s)$ in the following sense:
 $$\forall x\in\M, \quad x\in\rA
 \equival (\forall s\in\reel,\ \theta_s(x)=x)$$
Moreover $\M$ is semi-finite and can be canonically equipped with
a normal faithful semifinite trace $\tau$ such that
 $$\forall x\in\M,\quad \tau\rond \theta_s=e^{-s}\tau$$

Note that the von Neumann algebra $\M$ is independent from the
choice of the n. s. f. weight $w$ on $\rA$ (up to a
$*$-isomorphism preserving the trace and the group $(\theta_s)$).

  Let $L_0(\M,\tau)$ be the space of measurable operators associated with
$\tau$ (in
Nelson's sense [N]). Recall that
$L_0(\M,\tau)$ is the completion of $\M$, when $\M$ is equipped with the vector
space topology given by the neighborhoods of the origin:
 $$N(\eps,\delta)=\{x\in\M\mid \exists e\in\M
 \hbox{ projection s. t. }
 \N{xe}\leq \eps \hbox{ and }\tau(e^\perp)<\delta\}$$
Then the operations on $\M$ extend
by continuity to $L_0(\M,\tau)$, which becomes a topological *-algebra.

Note that if $\M$ acts on a Hilbert space $H$, $L_0(\M,\tau)$ can
be identified with a class of unbounded, closed, densely defined
operators on $H$ affiliated with $\M$. The operations on $L_0(\M,\tau)$
are identified with the strong sum and the strong product of unbounded
operators
(i. e. the sum, resp. the product followed by the closure operation).

If $h$ is an element of $L_0(\M,\tau)$, we define its {\it left
support} $\ell(h)$ (resp {\it right support} $r(h)$) as the least
projection $e$ of $\M$ such that $eh=h$ (resp. $he=h$). Clearly
$\ell(h^*)=r(h)$, so if $h$ is self-adjoint, $\ell(h)=r(h)$ which
we call  then simply the {\it support} of $h$ and denote by
$s(h)$.

The space $L_0(\M,\tau)$ is equipped with a positive cone
 $$L_0(\M,\tau)_+=\{h^*h\mid h\in L_0(\M,\tau)\}$$
which is the completion of the positive cone of $\M$. Every
element $h\in L_0(\M,\tau)$ has a unique polar decomposition
$$h=u\abs{h}$$
where $\abs{h}=(h^*h)^{1/2}\in L_0(\M,\tau)_+$ and $u$ is a
partial isometry of $\M$ whose right support is equal to that of
$h$.

The *-automorphisms $\theta_s$, $s\in\reel$  extend to
*-automorphisms of $L_0(\M,\tau)$. For $0<p\le\infty$, the space
$L_p(\rA)$ is defined by
 $$L_p(\rA)=\{h\in L_0(\M,\tau)\mid \theta_s(h)=e^{-s/p}h\}$$
The space $L_\infty(\rA)$ coincides with $\rA$ (modulo the
inclusions $\rA\subset\M\subset L_0(\M,\tau)$). The spaces
$L_p(\rA)$  are closed self-adjoint linear subspaces of
$L_0(\M,\tau)$. They are closed under left and right
multiplications by elements of $\rA$.  If $h=u\abs{h}$ is the
polar decomposition of $h\in L_0(\M,\tau)$, then
 $$h\in L_p(\rA)\equival u\in\rA\hbox{ and }\abs{h}\in L_p(\rA)$$
As a consequence, the left and right supports of $h\in L_p(\rA)$
belong to $\rA$.

It was shown by Haagerup that there is a linear homeomorphism
$\ph\mapsto h_\ph$ from $\rA_*$ onto $L_1(\rA)$ (equipped with the
vector space topology inherited from $L_0(\M,\tau) $), and this
homeomorphism preserves the additional structure (conjugation,
positivity, polar decomposition, action of $\rA$). It permits to
transfer the norm of $\rA_*$ to a norm on $L_1(\rA)$, denoted by
$\N{\ }_1$.

The space $L_1(\rA)$ is equipped with a distinguished bounded
positive linear form Tr, the ``trace'', defined by
 $$\forall\ph\in \rA_*,\quad {\rm Tr}\,(h_\ph)=\ph(\1)$$
Consequently, $\N{h}_1={\rm Tr}\,(\abs{h})$ for every $h\in
L_1(\rA)$.

For every $0<p<\infty$, the Mazur map $\rA_+\to\rA_+$, $x\mapsto
x^p$ extends by continuity to a map $L_0(\M,\tau)_+\to
L_0(\M,\tau)_+$, $h\mapsto h^p$. Then
 $$\forall h\in L_0(\M,\tau)_+\quad
 h\in L_p(\rA)\equival h^p\in L_1(\rA)$$
For $h\in L_p(\rA)$ set $\N{h}_p=\N{\,\abs{h}^p}_1^{1/p}$. Then
$\N{\ }_p$ is a norm when $1\leq p<\infty$, and a $p$-norm when
$0<p<1$ (see [Ko3] for this case). The associated vector space
topology coincides with that inherited from $L_0(\M,\tau)$.

Another important link between the spaces $L_p(\rA)$ is the {\it
external product}: in fact the product of $L_0(\M,\tau)$,
$(h,k)\mapsto h\cdot k$, restricts to a bounded bilinear map
$L_p(\rA)\times L_q(\rA)\to L_r(\rA)$, where $\frac 1r=\frac
1p+\frac 1q$. This bilinear map has norm one (``non commutative
H\"older inequality'').

Assume that $\frac 1p+\frac 1q=1$. Then the bilinear form
$L_p(\rA)\times L_q(\rA)\to \comp$, $(h,k)\mapsto {\rm
Tr}\,(h\cdot k)$ defines a duality bracket between $L_p(\rA)$ and
$L_q(\rA)$, for which $L_q(\rA)$ is (isometrically) the dual of
$L_p(\rA)$ (if $p\neq\infty$); moreover we have the tracial
property:
 $$\forall h\in L_p(\rA), k\in L_q(\rA),\quad
 {\rm Tr}\,(hk)={\rm Tr}\,(kh)$$

 \proclaim Definition 1.1. i) Two elements $h, k\in L_0(\M,\tau)$ are
called {\rm disjoint}, written as $h\perp k$, if they have
disjoint left, resp. right supports:\sn
 \centerline{$\ell(h)\perp \ell(k)\quad \hbox{and} \quad r(h)\perp  r(k)$}
\sn
  ii) A sequence $(h_n)\subset L_0(\M,\tau)$ is called {\rm
disjoint} if the $h_n$'s are pairwise disjoint; if in addition
$(h_n)\subset L_p(\rA)$ ($0<p<\infty$), $(h_n)$ is called {\rm
almost disjoint} if there is a disjoint sequence $(h_n')$ such
that $\lim_n\|h_n-h_n'\|_p=0$.

Note that if $(h_n)$ is almost disjoint in $L_p(\rA)$, so is
$(h_{\pi(n)})$ for every permutation $\pi$ on $\nat$; thus we can
speak of an almost disjoint countable subset in $L_p(\rA)$.

\smallskip

We shall repeatedly  use of the following two facts.

\proclaim Fact 1.2. i) If $h\in L_0(\M,\tau)_+$ and $0<p<\infty$,
then $s(h^p)=s(h)$.
 \hfill\break ii) If $h,k\in L_0(\M,\tau)$, then
$hk=0$ iff $r(h)\perp \ell(k)$.

\proof This is easy via a realization of $L_0(\M,\tau)$ as a set
of of unbounded, closed, densely defined operators on a Hilbert
space $H$. Then $\ell(h)$, resp. $r(h)^\perp$ is the projection
onto the closure of the range of $h$, resp. onto the kernel of
$h$. If $h\in L_0(\M,\tau)_+$, then it is self-adjoint and
property i) is well known. Concerning ii) we note that if $hk=0$,
then ${\rm ran}\,(k)\subset\ker h$, so $\ell(k)\leq r(h)^\perp$;
conversely if $r(h)\perp \ell(k)$ then $hk=hr(h)\ell(k)k=0$.\cqfd

\proclaim Fact 1.3. Let $0<p<\infty$ and $h,k$ be two elements of
$L_p(\rA)$.\hfill\break
 i) If $h\perp k$, then
$\N{h+k}_p^p=\N{h}_p^p+\N{k}_p^p$.\hfill\break
 ii) Conversely if
$p\neq 2$ and $\N{h+k}_p^p=\N{h-k}_p^p=\N{h}_p^p+\N{k}_p^p$,  then
$h\perp k$.

\proof i) If $h\perp k$, then $\abs{h+k}^p=\abs{h}^p+\abs{k}^p$,
hence ${\rm Tr}\abs{h+k}^p ={\rm Tr}\abs{h}^p+{\rm Tr}\abs{k}^p$.

ii) These two equalities implies that $h,k$ verify the equality
case of Clarkson's inequality. By Theorem A1 of the Appendix,
these elements are disjoint.\cqfd
\medskip
We finally mention the following ``localization'' fact which will
be used in section 2.

\proclaim Fact 1.4. If $e$ is an arbitrary projection of $\rA$,
then the subspace $eL_p(\rA)e$ is isometrically isomorphic to
$L_p(e\rA e)$, the $L_p$-space associated with the reduced von
Neumann algebra $e\rA e$; this isomorphism preserves the bimodule
structure (over $e\rA e$) as well as the external product in the
$L_p$ scale; in particular, it preserves the disjointness.

This is easily seen by taking a special n. s. f. weight of the form
$$w(x)=w_1(exe)+w_2(e^\perp x e^\perp)$$
where $w_1$, $w_2$ are n. s. f. weights on $e\rA e$, resp.
$e^\perp\rA e^\perp$. Then it is easy to see that $e$ is invariant
under the automorphism group $(\sigma_t^w)$ associated with $w$
and that $\sigma_t^{w_1}$ is nothing but the restriction of
$\sigma_t^w$ to $e\rA e$. Thus the crossed product $\M_e$
associated with $e\rA e$ is nothing but $e\M e$, on which the dual
automorphism group $(\theta_s^e)$ is simply the restriction of
$(\theta_s)$ and the trace $\tau_e$ is the restriction of $\tau$.
\goodbreak

\medskip \noindent{\bf Ultrapowers of non-commutative
$L_p$-spaces}
\medskip

Let $\U$ be an ultrafilter over some index set $I$.  If $X$ is a
Banach space, let $\ell_\infty(I;X)$ be the Banach space of
bounded families of elements of $X$, indexed by $I$, equipped with
the usual supremum norm. Let $N^\U$ be the subspace of
$\U$-vanishing families, i.e.
$$N^\U=\{(x_i)_{i\in I}\in\ell_\infty(I;X)\mid \lim_{i,\U}\N{x_i}=0\}$$
The ultrapower $X_\U$ is simply the quotient Banach space
$\ell_\infty(I;X)/N^\U$. If $(x_i)_{i\in I}$ is a member of
$\ell_\infty(I;X)$, we denote by $(x_i)^\bullet$ its image by the
canonical surjection $\ell_\infty(I;X)\to X_\U$. The norm of this
later element is simply given by $\Vert (x_i)^\bullet\Vert
=\lim\limits_{i,\U}\N{x_i}$.

The space $X$ is canonically isometrically embedded into its
ultrapower $X_\U$ by the diagonal embedding $x\mapsto
(x)_i^\bullet$, where $(x)_i$ is the constant family (all the
members of which are equal to $x$). We set $\hat x=(x)_i^\bullet$.
Sometimes we omit the hat over $x$ when no confusion can occur. If
$X$ is not finite dimensional and the ultrafilter is not trivial
(i.e. not principal) then $X\neq X_\U$.

If $X, Y$ are Banach spaces and $T: X\to Y$ is a bounded linear
operator, we can define canonically the ultrapower $T_\U$ of $T$
as the operator $X_\U\to Y_\U$, $(x_i)^\bullet\mapsto
(Tx_i)^\bullet$. More generally we can define analogously the
ultrapower $F_\U$ of a locally uniformly continuous map $F: X\to
Y$. In particular, if $B:X\times Y\to Z$ is a bounded bilinear
map, it has an ultrapower map $B_\U: X_\U\times Y_\U\to Z_\U$
defined by $B(\xi,\eta)=(B(x_i,y_i))^\bullet$ whenever
$\xi=(x_i)^\bullet, \eta= (y_i)^\bullet$.

All these are also valid for quasi-Banach spaces.

\smallskip

 Now let $\rA$ be
a C*-algebra. Then $\rA_\U$ is an involutive complex Banach
algebra when equipped with the natural product and conjugation
operations which are the respective ultrapowers of the product and
the conjugation operations of $\rA$. In fact, $\rA_\U$ is a
C*-algebra since it verifies the axiom $\N{xx^*}=\N{x}^2$ for
every $x\in\rA_\U$, which characterizes C*-algebras among
involutive complex Banach algebras. On the other hand, the class
of von Neumann algebras (dual C*-algebras) is not closed under
ultrapowers. However it was shown by U. Groh [G] that the class of
the preduals of von Neumann algebras is closed by ultrapowers. So
if $\rA$ is a von Neumann algebra, and $\rA_*$ is its (unique)
predual, then $(\rA_*)_\U$ is isometric to the predual of a von
Neumann algebra $\A$; moreover, $\rA_\U$ identifies naturally to a
w*-dense subspace  of $\A$. In fact, one can require that $\rA_\U$
be a *-subalgebra of $\A$; then $\A$ is uniquely defined as
C*-algebra. Note that the class of the preduals of semi-finite von
Neumann algebras is {\it not} closed under ultrapowers (see [Ra]),
which justifies the use of Haagerup $L_p$-spaces in the present
paper.

Groh's theorem was extended by the first author to the class of
non-commutative $L_p$-spaces, for arbitrary positive real $p$. It
was shown in [Ra] that $L_p(\rA)_\U$ is isometrically isomorphic
to $L_p(\A)$, where $\A$ does not depend on $p$ (it is precisely
the dual of $(\rA_*)_\U$). In fact, the isomorphisms
$\Lambda_p:L_p(\rA)_\U\to L_p(\A)$ constructed in [Ra] preserve
some more structures. Note that on $L_p(\rA)_\U$ there are
conjugation, absolute value map, left and right actions of
$\rA_\U$, and external products with other spaces $L_q(\rA)_\U$,
which are simply the ultrapowers of the corresponding operations
involving respectively $L_p(\rA)$, $\rA$ and $L_q(\rA)$.  Then the
identification maps $\Lambda_p$ preserve:

- conjugation: $\Lambda_p(\t h^*)=\Lambda_p(\t h)^*$

- absolute values: $\Lambda_p(\vert\t h\vert)=\vert\Lambda_p(\t
h)\vert$

- $\rA_\U$-bimodule structure: $\Lambda_p(\t x\cdot \tilde h\cdot
\t y)=\t x\cdot \Lambda_p(\tilde h)\cdot \t y$

- external product: $\Lambda_r(\t h\cdot \t k)=\Lambda_p(\t
h)\cdot \Lambda_p(\t k)$, $\frac 1r=\frac 1p+\frac 1q$

\noindent for all $\t h\in L_p(\rA)_\U$, $\t k\in L_q(\rA)_\U$,
$\t x,\t y\in \rA_\U$. We shall frequently use these properties
without any further reference.

\medskip \noindent {\bf $\ell_p$-sequences in Banach spaces}

\medskip
 A basic sequence $(x_n)$ of a Banach space $X$ is {\it
$K$-equivalent} to the unit $\ell_p$-basis iff there are positive
reals $a,b$ with $b/a\leq K$ such that
 $$a\,(\sum_n\abs{\lambda_n}^p)^{1/p}
 \leq\Vert\sum_n\lambda_nx_n\Vert
 \leq b\,(\sum_n\abs{\lambda_n}^p)^{1/p}$$
for every system $(\lambda_n)$ of finitely nonzero complex
numbers.

It is {\it almost 1-equivalent} to the $\ell_p$-basis if for some
sequence $(\eps_n)$ of positive reals such that
$\lim\limits_{n\to\infty}\eps_n=0$, the tail $(x_m)_{m\geq n}$ is
$(1+\eps_n)$-equivalent to the $\ell_p$-basis.

It is {\it asymptotically 1-equivalent} to the $\ell_p$-basis if for
some sequence $(\eps_n)$ of positive reals such that
$\lim\limits_{n\to\infty}\eps_n=0$ we have:
 $$(\sum_n(1-\eps_n)^p\abs{\lambda_n}^p)^{1/p}
 \leq\Vert\sum_n\lambda_n x_n\Vert
 \leq\,(\sum_n(1+\eps_n)^p\abs{\lambda_n}^p)^{1/p}$$
for every system $(\lambda_n)$ of finitely nonzero complex
numbers. The space spanned by such a sequence is called an asymptotically
isometric
copy of $\ell_p$ in the terminology of [DJLT]. Note that the subspace
spanned by a
sequence which is almost 1-equivalent to the $\ell^p$ does not contain
necessarily
an  asymptotically isometric copy of $\ell_p$ (see [DJLT]).

\bigskip\centerline{\bf 2. Elements in $L_p(\rA)_\U$ which are
disjoint from $L_p(\rA)$}
\medskip

 In this section we develop the main tools of the paper (Theorem
 2.3, Lemma 2.6 and Theorem 2.7). We also give several
 characterizations of bounded sequences in $L_p(\rA)$ which
 have almost disjoint subsequences.

\medskip\noindent {\bf Pairs of disjoint elements in ultrapowers}
\medskip

 Let $\rA$ be a von Neumann algebra, $\rA_*$ its predual
and $(\rA_*)_\U$ the ultrapower of $\rA_*$ relative to an
ultrafilter over some index set $I$. Let $\A$ be the dual von
Neumann algebra of $(\rA_*)_\U$. In this section we shall prove
that disjoint elements of $\A_*$, when considered in $(\rA_*)_\U$,
admit pairwise disjoint families of representatives in
$\ell_\infty(I;\rA_*)$. This property is easy to prove in the
commutative case, using the lattice operations in $L_1$-spaces.
The proof we give in the noncommutative case is based on the fact
that elements of the algebra $\A$ can be ``locally'' identified
with elements of the ultrapower $\rA_\U$.

Recall that a projection $p$ in a von Neumann algebra $\M$ is
{\it$\sigma$-finite} if it is the support of a normal state.
Equivalently, there is $h\in L_2(\M)_+$ such that $s(h)=p$.

\proclaim Proposition {2.1}. For every $x\in\A$ and every
$\sigma$-finite projection $p$ of $\A$ there is a family
$(x_i)\subset\rA$ with $\N{x_i}\leq\N{x}$ for every $i\in I$ and
representing an element $\t x=(x_i)^\bullet$ of $\rA_\U$ such that
$\t xp=xp$.

\proof  Before to start the proof, recall that the positive cone
in $L_1(\rA)_\U$ consists of elements representable by a bounded
family of nonnegative elements of $L_1(\rA)$.

\noindent Let now $\t k\in L_2(\A)=L_2(\rA)_\U$, $\t k\geq 0$, with
support $p$, and set $\t h= x\t k$. Note that
 $$0\leq \t h^*\t h=\t kx^*x\t k\leq \N{x}^2\t k^2$$
Let $(h_i)_{i\in I}$ be a bounded family in $L_2(\rA)$
representing $\t h$. We can find a bounded family $(\ell_i)_{i\in
I}$ in $L_1(\rA)_+$, representing $\t k^2$ and such that
 $$\forall i\in I,\quad h_i^*h_i\leq \N{x}^2\ell_i$$
For, let $(a_i)_{i\in I}$ be a bounded family in $L_1(\rA)$
representing $\ds\t k^2-{\t h^*\t h\over \N{x}^2}$: since this
later element is positive, we can choose $a_i\geq 0$ for every
$i\in I$; then set $\ds\ell_i=a_i+{h_i^*h_i\over \N{x}^2}$.

\noindent Then for every $i\in I$ there exists $x_i\in \rA$ such
that
 $$h_i=x_i\ell_i^{1/2}\hbox{ and } \N{x_i}\leq\N{x}$$
The bounded family $(\ell_i^{1/2})_{i\in I}$ represents $(\t
k^2)^{1/2}=\t k$, and so $x\t k=\t h=(x_i\ell_i^{1/2})^\bullet=\t
x\t k$, which implies $\t xp=xp$.\cqfd

\proclaim Corollary {2.2}. For every $x\in\A$, $x\geq 0$ and every
$\sigma$-finite projection $p$ of $\A$ there exists a family
$(x_i)_{i\in I}\subset\rA$ with $0\leq x_i\leq \N{x}$ representing
an element $\t x$ of $\rA_\U$ such that $p\t xp=pxp$.

\proof Applying  Proposition 2.1 to $y=x^{1/2}$ and $p$, we obtain
$(y_i)\subset\rA$, with $\N{y_i}\leq\N{y}=\N{x}^{1/2}$ and $\t
y=(y_i)^\bullet$ satisfying $\t y p=y p$. Set $x_i=y_i^*y_i$, then
$p\t x p=p\t y^*\t y p=p y^2p=pxp$. \cqfd

\smallskip
 The next result states that two disjoint $\sigma$-finite
projections of $\A$ can be separated by a projection of $\rA_\U$.
It is the key technical result of the paper.

\proclaim Theorem {2.3}. Let $p,q$ be two disjoint $\sigma$-finite
projections in $\A$. There exists a family of projections
$(r_i)_{i\in I}$ in $\rA$ representing a projection $\t r$ of
$\rA_\U$ such that:\smallskip \centerline{$\t r\geq p$~ and~~ $\t
r^\perp\geq q$}

\proof Applying the preceding corollary to $x=p$ and the
$\sigma$-finite projection $s=p+q$, we find $(x_i)_{i\in
I}\subset\rA$ with $0\leq x_i\leq  \1$ such that $\t
x=(x_i)^\bullet$ verifies:
 $$s\t x s=sps=p  $$
Then
 $$s(\1-\t x) s=s-p=q$$
Hence:
 $$\left\{\eqalign{\t x s&=p+s^\perp\t x s\cr
 (\1-\t x)s&=q+s^\perp(\1-\t x) s}\right.$$
whence:
 $$\left\{\eqalign{\t x p&=p+s^\perp\t x p\cr
 (\1-\t x)q&=q+s^\perp(\1-\t x) q}\right.$$
Let $\t k\in L_2(\A)=L_2(\rA)_\U$ with support $p$. Since $\N{\t
x}\leq 1$, we have
 $$\eqalign{\Vert\t k\Vert_2^2\geq \Vert\t xs \t k\Vert^2_2
 &=\Vert p\t k\Vert_2^2+\Vert s^\perp\t x p\t k\Vert_2^2\cr
 &= \Vert \t k\Vert_2^2+\Vert s^\perp\t x \t k\Vert_2^2}$$
Therefore, $s^\perp \t x \t k=0$, and so $s^\perp \t x p=0$.
Similarly, since $\N{\1-\t x}\leq 1$, we have  $s^\perp(\1- \t x)
q=0$. Thus:
 $$\t x p=p~\hbox { and }  (\1-\t x)q=q$$
For every $i\in I$ consider the spectral projection
$r_i=\chi_{[\frac 12,1]}(x_i)$ of $x_i$ associated with the
indicator function of the interval $[\frac 12,1]$. Note that
$r_i=f(x_i)x_i$, where the function $f$ is defined by $\ds
f(t)=t^{-1} \chi_{[\frac 12,1]}(t)$. We have
$\N{f(x_i)}\leq\N{f}_\infty=2$, so the family $(f(x_i))_{i\in I}$
defines an element of $\rA_\U$. Then:
 $$(r_i)^\bullet q=(f(x_i))^\bullet(x_i)^\bullet q=0$$
Similarly, since $r_i^\perp=g(\1-x_i)(\1-x_i)$ with $\ds g(t)=t^{-1}
\chi_{(\frac 12,1]}(t)$, we deduce
 $$(r_i^\perp)^\bullet p=(g(\1-x_i))^\bullet (\1-x_i)^\bullet p=0$$
Therefore $\t r=(r_i)^\bullet$ is a desired projection of
$\rA_\U$. \cqfd

\smallskip

\proclaim Corollary {2.4}. Let $0<p<\infty$. Two elements $\t h,\t
k$ of $L_p(\A)=L_p(\rA)_\U$ are disjoint if and only if they admit
representative families $(h_i)_{i\in I}$, $(k_i)_{i\in I}$ such
that for every $i\in I$, $h_i$ is disjoint from $k_i$

\proof The ``if'' part is evident. To prove the necessity of the
condition assume $\t h,\t k$ are disjoint. By Theorem 2.3, we can
find projections $\t r=(r_i)^\bullet, \t s=(s_i)^\bullet$ in
$\rA_\U$ such that
 $$\eqalign{
 &\t s\geq \ell(\t h),~~~~~~~~~\t s^\perp\geq \ell(\t k)\cr
 & \t r\geq r(\t h),~~~~~~~~~ \t r^\perp\geq r(\t k)}$$
Let $(h_i)_{i\in I}$, resp $(k_i)_{i\in I}$ be two representative
families of $\t h$, resp $\t k$. Let $h'_i=s_ih_ir_i$ and
$k'_i=s_i^\perp k_ir_i^\perp$. Then $(h'_i), (k'_i)$ are two
desired representative families.\cqfd

\medskip \noindent {\bf Elements of $L_p(\rA)_\U$ disjoint from
$L_p(\rA)$} \medskip

We say that $\t h\in L_p(\rA)_\U$ is disjoint from $L_p(\rA)$
(considered as subspace of $L_p(\rA)_\U$) if $\t h\perp k$ for
every $k\in L_p(\rA)$. (This is an abuse of notation; we should
write $\t h\perp \hat k$, where $\hat k=(k)^\bullet$ is the
canonical image of $k$ in $L_p(\rA)_\U$; note that the left and
right supports of $\hat k$ do not coincide with the canonical
images in $\rA_\U$ of the supports of $k$). Equivalently, $\t h
k=k\t h=0$ for every $k\in L_p(\rA)$. Since $k\in L_p(\rA)_+$ and
$k^\alpha\in L_{p/\alpha}(\rA)_+$ have the same support (in $\A$),
another equivalent condition is that $k\t h=0=\t h k$ for every
$k$ in $L_q(\rA)$, for some (every) $q$, $0<q<\infty$.

A simple example is given by the following lemma:

\proclaim Lemma {2.5}. Suppose that $\U$ is a free ultrafilter
over $\nat$ and let $(h_n)_{n\in\nat}$ be a bounded disjoint
sequence in $L_p(\rA)$. Then the element $\t h$ defined by this
sequence in $L_p(\rA)_\U$ is disjoint from $L_p(\rA)$.

\proof The left supports $s_n=\ell(h_n)$ are pairwise disjoint. If
$k\in L_2(\rA)$ we have then $\N{ks_n}_2\to 0$. For, since the
elements $ks_n$, $n\in\nat$ are pairwise orthogonal for the
natural scalar product of $L_2(\rA)$:
 $$\sum_n\N{ks_n}_2^2=\Vert\sum_n ks_n\Vert_2^2\leq \N{k}_2^2$$
Consequently, by the H\"older inequality (with $1/r=1/2+1/p$),
 $\|kh_n\|_r\leq\N{ks_n}_2\N{h_n}_p\to 0$. Similarly,
$\N{h_nk}_{r}\to 0$. A fortiori,
$\lim\limits_{n,\U}\N{h_nk}_{r}=0=
\lim\limits_{n,\U}\N{kh_n}_{r}=0$, which implies  $k\t h=\t hk=0$.
\cqfd

\smallskip

\proclaim Lemma {2.6}. Let $0<p<\infty$, and let $\cal S$ be a
separable subset of elements of $ L_p(\rA)_\U$ which are disjoint
from $L_p(\rA)$. For each $\t h\in \cal S$ let  $(h_i)_{i\in I}$
be a bounded family in $ L_p(\rA)$ defining $\t h$.
 Then for every finite system $\P$ of pairwise commuting projections of $\rA$
and every separable subset $\K$ of $L_p(\rA)$ there exists a
family $(s_i)$ of projections of $\rA$ commuting with $\P$ and
such that:\smallskip \centerline{$\left\{\matrix{&\tvi\forall
k\in\K,\hfill&\qquad\|s_ik\|_p+\N{ks_i}_p\To\limits_{i,\U}0\hfill\cr
&\forall \t h\in{\cal S},&\qquad\N{s_i^\perp
h_i}_p+\N{h_is_i^\perp }_p\To\limits_{i,\U}0\hfill\cr}\right.$}

\medskip \proof  Let $\P=\{p_1,...,p_N\}$: replacing
$\P$ by the set of atoms of the (finite) Boolean algebra generated
by $\P$, we may suppose that the $p_j$ 's are disjoint and
$\sum\limits_{j=1}^N p_j=\1$ . Note that for every $j=1,...,N$, and
$\t h\in \cal S$ the elements $\hat p_j\t h$ and $\t h \hat p_j$
are disjoint from $L_p(\rA)$, and a fortiori from $\hat p_j
L_p(\rA)\hat p_j$. We may identify $p_jL_p(\rA)p_j$ with
$L_p(p_j\rA p_j)$, and $\hat p_jL_p(\rA)_\U\hat p_j$ with
$L_p(p_j\rA p_j)_\U$.  Let
 $$e=\bigvee\limits_{\t h\in{\cal S}}\ell(\hat p_j\t h)\vee
 r(\t h\hat p_j)\qquad
 f=\bigvee\limits_{k\in\K}\ell(\hat p_j\hat k)\vee
 r(\hat k\hat p_j)$$
Then $e$ and $f$ are $\sigma$-finite disjoint projections. Note
that all the support projections above are smaller than $\hat
p_j$, hence belong to $\hat p_j\A\hat p_j$, and so do $e$ and $f$.
Thus by Theorem 2.3 there exists a family $(s_i^{(j)})_i$ of
projections of $p_j\rA p_j$ such that the corresponding
projections $\t s^{(j)}$ of $\hat p_j\A\hat p_j$ satisfy $e\le \t
s^{(j)}$ and $f\le \big(\t s^{(j)}\big)^\perp$.

We set
 $$s_i=\sum_{j=1}^Ns^{(j)}_i= \sum_{j=1}^Np_js^{(j)}_ip_j.$$
Then all $s_n$ clearly commute with $\P$, and
 $$\eqalign{(s_ih_i)^\bullet&=\sum_{j=1}^N\t s^{(j)}\t h
 =\sum_{j=1}^N\t s^{(j)}\hat p_j\t h
 =\sum_{j=1}^N\hat p_j\t h=\t h,
 \hbox{ for every } \t h\in {\cal S}\cr
 (s_ik)^\bullet&=\sum_{j=1}^N\t s^{(j)}\t k
 = \sum_{j=1}^N\t s^{(j)}\hat p_j\t k=0,
 \hbox{ for every } k\in \K}$$
Similarly, $(h_is_i)^\bullet=\t h$ and $(ks_i)^\bullet=0$ for all
$\t h\in {\cal S}$ and  $k\in\K$. Therefore, the family $(s_i)$
satisfies all requirements of the lemma. \cqfd

\smallskip
 Recall  that a sequence $(h_n)$ in $L_p(\rA)$ is almost disjoint
 if there is a disjoint sequence $(h_n')\subset L_p(\rA)$ such that
 $\lim_n\|h_n-h_n'\|_p=0$ (see Definition 1.1).

\proclaim Theorem {2.7}. A bounded family $(h_i)_{i\in I}$ in
$L_p(\rA)$ has an almost disjoint countable subfamily if and only
if for some free ultrafilter $\U$ over $I$ $(h_i)_{i\in I}$
defines an element of the ultrapower $L_p(\rA)_\U$ which is
disjoint from $L_p(\rA)$.

\proof The ``only if'' part results from  Lemma 2.5, choosing an
ultrafilter $\U$ containing as an element the infinite subset of $I$
indexing the countable subfamily. Let us prove the ``if'' part.

We use Lemma 2.6 to construct inductively a sequence of distinct
indices $(i_n)$ and a sequence $(q_n)$ of commuting projections of
$\rA$, such that
 $$\forall n\in\nat: \Vert q_n^\perp h_{i_n}\Vert_p
 +\Vert h_{i_n}q_n^\perp\Vert_p<2^{-n}\hbox{ and }\forall m< n,\
 \N{q_nh_{i_m}}_p+\N{h_{i_m}q_n}_p<2^{-n}$$
We start with some $i_0\in I$ and $q_0=\1$. At the $(n+1)$-th step
apply Lemma 2.6 to $\P=\{q_0,...,q_n\}$ and
$\K=\{h_{i_0},...,h_{i_n}\}$.

Set $p_n=q_n(\bigwedge\limits_{k>n}q_k^\perp)$. The projections
$p_n$ are disjoint. Note that since the $q_k$'s commute, we have
$\bigvee\limits_{k>n}q_k=\sum\limits_{k>n}x_nq_k$ for some
$x_n\in \rA$, $0\leq x_n\leq \1$. Then:
 $$\N{(p_n^\perp-q^\perp_n)h_{i_n}}_p=\N{(q_n-p_n)h_{i_n}}_p
 =\Vert q_n(\bigvee_{k>n} q_k)h_{i_n})\Vert_p
 \leq \Vert (\sum_{k>n}x_nq_k)h_{i_n} \Vert_p$$
 Therefore, if
$p\geq 1$,
 $$\N{(p_n^\perp-q^\perp_n)h_{i_n}}_p\le
  \sum\limits_{k>n}\Vert q_k h_{i_n}\Vert_p\le 2^{-n}\,;$$
similarly, if $0<p<1$,
 $$\N{(p_n^\perp-q^\perp_n)h_{i_n}}_p
 \le(2^p-1)^{-1/p}2^{-n}.$$
Thus in both cases,
$\N{(p_n^\perp-q^\perp_n)h_{i_n}}_p\To\limits_{n\to \infty} 0$.
Hence it follows that $\N{p_n^\perp h_{i_n}}_p\to 0$. In the same
way, we show that $\N{h_{i_n}p_n^\perp}_p\to 0$. Therefore,
$\N{h_{i_n}-p_n h_{i_n}p_n}_p\to 0$.\cqfd

\smallskip
\noindent {\bf Remarks {2.8}:} i) Theorem 2.7 has a close analog
for left (resp. right) disjointness. Say that a sequence $(h_n)$
in $L_p(\rA)$ is almost left (resp. right) disjoint if there
exists a sequence $(h'_n)$ of pairwise left (resp. right) disjoint
vectors in $L_p(\rA)$ such that $\N{h_n-h'_n}_p\to 0$. Similarly,
an element $\t h$ of the ultrapower $L_p(\rA)_\U$ is left (resp.
right) disjoint from $L_p(\rA)$ if it is left (resp. right)
disjoint from every element of $L_p(\rA)$ (canonically embedded
in) $L_p(\rA)_\U$. Then a bounded family in $L_p(\rA)$ has an
almost left (resp. right) disjoint countable subfamily iff for
some free ultrafilter $\U$ over the index set $I$ it defines an
element of the ultrapower $L_p(\rA)_\U$ which is left (resp.
right) disjoint from $L_p(\rA)$.

ii) The proof of Theorem 2.7 shows in fact the following: if
$(h_n)$ is an almost disjoint sequence in $L_p(\rA)$, then there
exist a subsequence $(h_{i_n})$ and a disjoint sequence $(p_n)$ of
projections in $\rA$ such that $\N{h_{i_n}-p_n h_{i_n}p_n}_p\to
0$. The same remark also applies to left and right almost disjoint
sequences.

\medskip\noindent {\bf Disjoint types over $L_p(\rA)$}

\medskip

 Recall that following [KrM] a {\it type} over a Banach
space $E$ (or $p$-Banach space) is a function $\tau: E\to \reel_+$
of the form $\tau(x)=\lim\limits_{i,\U}\N{x+x_i}$, where $(x_i)$
is a bounded family of points of $E$ and $\U$ an ultrafilter over
$I$. Equivalently,  we have $\tau(x)=\N{x+\xi}$, where $\xi$ is an
element of the ultrapower $E_\U$: then we say that $\xi$ defines
the type $\tau$. We say that {\it a sequence $(x_n)\subset E$
defines the type $\tau$ if
$\tau(x)=\lim\limits_{n\to\infty}\N{x+x_n}$}. Note that in
separable spaces every type is definable by a sequence. We call a
type $\tau$ over $L_p(\rA)$ a {\it disjoint type} if it is
definable by an element $\xi$ of some ultrapower $L_p(\rA)_\U$
which is disjoint from $L_p(\rA)$.

\proclaim Proposition {2.9}. If $0<p<\infty$,  the disjoint types
over $L_p(\rA)$ are exactly the functions of the form $h\mapsto
F_a(h)=(\N{h}^p+a^p)^{1/p}$, where $a$ is a nonnegative real
number. Moreover if $p\neq 2$ an element $\xi$ of an ultrapower
$L_p(\rA)_\U$ defines a disjoint type over $L_p(\rA)$ if and only
if it is disjoint from $L_p(\rA)$.

\proof If $\tau$ is a disjoint type defined by $\xi\in
L_p(\rA)_\U$ with $\xi\perp L_p(\rA)$, we clearly have  $\tau=F_a$
with $a=\N{\xi}$. Conversely, if $\xi$ defines the type $F_a$,
then $a=\N{\xi}^p$ and for every $h\in L_p(\rA)$:
$$\N{h\pm\xi}^p=F_a(\pm h)^p=(\N{h}^p+\N{\xi}^p)$$
hence $\xi\perp h$ by Fact 1.3 when $p\neq 2$. \cqfd

\proclaim Lemma 2.10. Every disjoint normalized sequence $(h_n)$
in $L_p(\rA)$ defines a disjoint type.

\proof Let $\U$ be a free ultrafilter over $\nat$. By lemma 2.5, $\t
h=(h_n)^\bullet$ is disjoint from $L_p(\rA)$. Hence for every
$k\in L_p(\rA)$, $\lim\limits_{n,\U}\N{k+h_n}=\Vert k+\t
h\Vert=(\N{k}^p+\Vert\t h\Vert^p)^{1/p}=F_1(k)$. Since this is
true for every ultrafilter $\U$, we have
$\lim\limits_{n\to\infty}\N{k+h_n}=F_1(k)$.\cqfd

The following gives several characterizations of a bounded
sequence which defines a disjoint type:

\proclaim Proposition {2.11}. Let $0<p,q<\infty$, $p\neq 2$ and
$(h_n)$ be a bounded sequence in $L_p(\rA)$. Assume that the
sequence of norms $(\Vert h_n\Vert)$ converges. Then the following
assertions are equivalent:\hfill\break
 i) $(h_n)$ defines a disjoint type.\hfill\break
 ii) For every element $h$ of $L_q(\rA)$
we have $\lim\limits_{n\to\infty}h\cdot h_n=0=
\lim\limits_{n\to\infty}h_n\cdot h$.\hfill\break
 iii) Every subsequence of
$(h_n)$ contains a subsequence which is almost disjoint in
$L_p(\rA)$.\hfill\break
 iv) Every subsequence of $(h_n)$ contains a
subsequence asymptotically 1-equivalent to the $\ell_p$-basis (up
to a constant factor).\hfill\break
 v) (For $p\geq 1$:) Every
subsequence of $(h_n)$ contains a subsequence asymptotically
1-equivalent to the $\ell_p$-basis (up to a constant factor) and
spanning an almost complemented subspace of $L_p(\rA)$.

\proof Note that the hypothesis on the convergence of the norms is
necessary since if $(h_n)$ defines a disjoint type $F_a$, then
$\N{h_n}\to F_a(0)=a$.

\noindent i) $\implies$ ii): For every free ultrafilter $\U$ over
$\nat$, the element $\t h$ defined by $(h_n)$ in $L_p(\rA)_\U$
defines the same disjoint type. By Proposition 2.9, $\t h$ is
disjoint from $L_p(\rA)$. Equivalently, for every $k\in L_q(\rA)$
we have $\hat k\t h=0=\t h\hat k$, where $\hat k=(k)^\bullet$ is
the canonical image of $k$ in $L_q(\rA)_\U$. Since $\hat k \t
h=(kk_n)^\bullet$, $\t h\hat k=(h_nk)^\bullet$, we have
$\lim\limits_{n,\U} kh_n=0 =\lim\limits_{n,\U} h_nk$. Since this
is true for every free ultrafilter $\U$ over $\nat$, we have
$\lim\limits_{n\to\infty} kh_n=0 =\lim\limits_{n\to\infty} h_nk$.

\noindent ii) $\implies$ iii): Since every subsequence of $(h_n)$
verifies also hypothesis ii), we may argue with the whole
sequence. Let $\U$ be a free ultrafilter over $\nat$  and $\t h$
the element defined by $(h_n)$ in $L_p(\rA)_\U$. Then $\t h$ is
disjoint from $L_p(\rA)$, and by Theorem 2.7 $(h_n)$ has an almost
disjoint subsequence.

\noindent iii) $\implies$ iv) (resp. and v) if $p\geq 1$): It is
clear that every sequence of normalized pairwise disjoint elements of
$L_p(\rA)$ is isometrically equivalent to the $\ell_p$-basis
(resp. and spanning a $1$-complemented subspace if $p\geq 1$). By
standard perturbation techniques (see e. g. [LT], prop. 1. a. 9)
one deduces that an almost disjoint sequence of elements whose
norm converges to $1$ has a subsequence which is almost $1$-equivalent to the
$\ell_p$-basis (resp. and spanning an almost $1$-complemented
subspace). This subsequence $(h'_n)$ in turn has a subsequence which is
asymptotically
$1$-equivalent to the $\ell_p$-basis: This is a consequence of the fact that
$$\forall h\in L_p(\rA),\quad\lim_{n\to\infty}\N{h+h'_n}=(\N{h}^p+1)^{1/p}$$
and by a standard Ascoli type argument, this limit is uniform on the unit
ball of
every finite dimensional subspace $V$ of $L_p(\rA)$. Hence for every
$\delta>0$ there
exists $N=N(V,\delta)$ such that
$$\forall n\geq N,\forall h\in V,\forall
\lambda\in\comp, \quad (1-\delta)(\N{h}^p+\abs{\lambda}^p)\leq
\N{h+\lambda h'_n}^p\leq (1+\delta)(\N{h}^p+\abs{\lambda}^p)$$
Choose a sequence $(\delta_n)$ with $0<\delta_n<1$ and
$\prod\limits_n(1-\delta_n)>0$, and define inductively
$n_0=1<n_1<..<n_k<n_{k+1}...$
by applying the preceding to $V_k=\span[h_{n_k}\mid 1\leq \ell\leq k\,]$
and setting
$n_{k+1}=\max\,[n_k+1,N(V_k,\delta_k))]$.

\noindent iv) $\implies$ iii): Suppose that $h_n$ itself is
asymptotically equivalent to the $\ell_p$-basis. Let $\U$ be a
free ultrafilter over $\nat$ and for every $m\in \nat$ let $\t
h_m$ be the element of $L_p(\rA)_\U$ defined by the sequence
$(h_{m+n})_{n\in\nat}$. Then the sequence $(\t h_m)$ is
isometrically equivalent to the $\ell_p$-basis. Let us make the
identification $L_p(\rA)_\U=L_p(\A)$. By Fact 1.3 the elements $\t
h_m$ are disjoint in $L_p(\A)$. Let $\xi$ be the element of
$L_p(\A)_\U$ defined by the sequence $(\t h_m)$. It is disjoint
from $L_p(\A)$, hence a fortiori from $L_p(\rA)$. On the other
hand, iterated ultrapowers are ultrapowers (relative to the
product ultrafilter), i.e. $\xi\in L_p(\rA)_{\U\times\U}$. Recall
that the ultrafilter $\U\times\U$ over $\nat\times\nat$ is defined
by
 $$A\in\U\times\U\Ssi
 \{n\mid\{m\in\nat\mid (n,m)\in A\}\in\U\}\in\U$$
With this identification we have $\xi=(h_{m+n})_{m,n}^\bullet$.
Let $f: \nat\times\nat\to\nat$  be a bijective map and $\V$ be the
ultrafilter $f(\U\times\U)=\{f(A) \mid A\in\U\times\U\}$. Set
$\ph(f(m,n))=m+n$. Then the sequence $(h_{\ph(i)})$ defines in
$L_p(\rA)_\V$ an element disjoint from $L_p(\rA)$, so by Theorem
2.7 it has an almost disjoint subsequence. Since clearly
$\ph(i)\to\infty$ when $i\to\infty$, a further subsequence will be
a subsequence of the initial sequence $(h_n)$.

\noindent iii)$\implies$ i): by Lemma 2.10. \cqfd
\bigskip
\noindent{\bf  Remark.} In the case $p=2$, the relations
ii) $\equival$ iii) $\implies$ i) $\equival$ iv) $\equival$ v) are still
true.

\bigskip\centerline{\bf 3. Embedding of $\ell_p$-sums of finite
dimensional spaces}\medskip

 We begin this section by recalling some standard notions from
 Banach space theory.
If $(F_n)$ is a (finite or infinite) sequence of Banach (or
quasi-Banach) spaces, their $p$-direct sum $(\bigoplus_{n}^{}
F_n)_p$ is the space of sequences $(x_n)\in\prod\limits_n F_n$
such that $\sum\limits_n \N{x_n}_{F_n}^p$ converges, equipped with
the natural norm $\N{(x_n)}=(\sum\limits_n
\N{x_n}_{F_n}^p)^{1/p}$. As usual, if all the spaces $F_n$
coincide with a given space $F$, their $p$-direct sum is denoted
by $\ell_p(F)$ for an infinite sequence, $\ell_p^k(F)$ for a
finite sequence with $k$ elements.

A  Banach (or quasi-Banach) space  $X$ contains uniformly a
sequence of Banach (or quasi-Banach) spaces  $(Y_n)$ if for some
constant $K$ the space $X$ contains for every $n$ a subspace
$X_n$ which is $K$-isomorphic to $Y_n$; then we say that $X$
contains the $Y_n$'s $K$-uniformly.

We say that a sequence $(E_n)$ of closed subspaces of the space
$L_p(\rA)$ is {\it almost disjoint} if there exists a sequence
$(p_n)$ of pairwise disjoint projections of $\rA$ such that
 $\lim\limits_{n\to\infty}\N{T_n}=0$, where
 $T_n$ is the operator $E_n\to L_p(\rA)$, $x\mapsto
(x-p_nxp_n)$.

\smallskip

 The following is one of the main results of this paper.

\proclaim Theorem 3.1. Let $0<p<\infty$, $p\neq 2$, $\rA$ be a von
Neumann algebra and $X$  a closed subspace of $L_p(\rA)$. Let
$(F_n)$ be a sequence of finite dimensional normed (or
quasi-normed) spaces. \hfill\break
 i) If $X$ contains
$K$-uniformly the finite $p$-direct sums $\ell_p^n(F_j)$, $n,j\geq
1$, then $X$ contains a subspace isomorphic to the infinite
$p$-direct sum $(\bigoplus_{j}^{} F_j)_p$. \hfill\break
 ii) More
precisely, under the assumption of i) for every $\eps>0$ there
exists an almost disjoint sequence $(E_n)$ of finite dimensional
subspaces of $X$ such that for every $n$, $E_n$ is
$(K+\eps)$-isomorphic to $F_n$. \hfill\break
 iii) If in addition
$1\leq p<\infty$ and $X$ contains the $\ell_p^n(F_j)$ ($n,j\geq
1$) as uniformly complemented subspaces of $L_p(\rA)$, then the
$E_n$ can be found uniformly complemented, and consequently $X$
contains $(\bigoplus_{j}^{} F_j)_p$ as  a complemented subspace of
$L_p(\rA)$.

This result immediately implies Theorem 0.1:\medskip

\noindent{\bf Proof of Theorem 0.1:} To deduce Theorem 0.1 from
Theorem 3.1 we need only to note that for every $n,m\geq 1$, the
space $S_p^{nm}$ contains isometrically $\ell_p^n(S_p^m)$ (as
$1$-complemented subspace when $p\geq 1$) by just taking the
block-diagonal embedding.\cqfd

\smallskip

The remainder of this section is devoted to the proof of Theorem
3.1. We first refine the given embeddings of $\ell_p^n(F_j)$ into
$X$. We denote by $(e_i)$ the natural basis of $\ell_p$ or
$\ell_p^n$. If $F$ is a space and $x\in F$, then $e_i\otimes x$
denotes the sequence $(0,0,..., 0,x,0...)$, where $x$ is at the
$i$-th place.

\proclaim Lemma 3.2. Let $X$ be a ($p$-)Banach space and $F$  a
finite dimensional (quasi-)normed space. Assume that $X$ contains
$K$-uniformly the spaces $\ell_p^n(F)$, $n\geq 1$.  Then for every
$\eps>0$ and every $n\geq 1$  there exists a $K$-isomorphic
embedding $T_{n,\eps}: \ell_p^n(F)\hookrightarrow X$ such that for
every nonzero $x\in F$ the sequence
 $\ds\big(\N{T_{n,\eps}(e_i\otimes x)}^{-1}
 T_{n,\eps}(e_i\otimes x)\big)_{1\leq i\leq n}$
is $(1+\eps)$-equivalent to the unit basis of $\ell_p^n$. If in
addition $p\geq 1$ and the initial copies of the $\ell_p^n(F)$ are
$C$-complemented, the new ones $T_{n,\eps}(\ell_p(F))$ are
$KC$-complemented.

\proof Given $n\geq 1$  let $T_n$ be a $K$-isomorphic embedding of
$\ell^n_p(F)$ into $X$. We define canonically a $K$-isomorphic
embedding $T$ of $\ell_p(F)$ into some ultrapower $X_\U$ of $X$ by
extending the $T_n$ to operators $ \ell_p(F)\to X$ (simply set
$T_n(e_k\otimes x)=0$ if $k>n$) and then setting $T(e_k\otimes
x)=(T_n(e_k\otimes x))^\bullet$.

If $p\geq 1$, we use Krivine's Theorem (see [MS], Theorem 12.4 in the
real case; [BL], Ch. 6, Cor. 3 for the complex version): every basic
sequence
$(x_n)$ in a Banach space  contains, for some
$q\in [1,+\infty]$, every $\eps>0$ and every $n\geq 1$,  a finite
sequence of disjoint blocks  $(1+\eps)$-equivalent  to the
$\ell_q^n$-basis. Of course, if the sequence $(x_n)$ is itself
$K$-equivalent to the $\ell_p$-basis, then $q=p$.

If $0<p\leq 1$, we use the the following $p$-normed space version
of the well known  James distorsion theorem on $\ell_1$: if a
$p$-Banach space has a basis equivalent to the $\ell_p$-basis,  it
contains for every $\eps>0$ a basic sequence which is
$(1+\eps)$-equivalent to the $\ell_p$-basis, and consists of
successive blocks of the initial basis. We refer to [J] or to
[LT], Proposition 2e3 for the proof of James' Theorem in the case
$p=1$. The adaptation of this proof to the case $0<p<1$ is
straightforward.

Fix a non-zero $\xi\in F$. In both cases, for every $n\geq 1$ we
can find a sequence $J_1<J_2<...<J_n$ of successive disjoint
intervals of $\nat$ and systems $(\lambda_{k,j})_{k=1,...,n; j\in
J_k}$ such that $\sum_{j\in J_k}\vert\lambda_{k,j}\vert^p=1$ for
every $k$, and such that
 $$\forall\rho_1,...,\rho_n\in\comp, \quad
 \Vert\sum_{k=1}^n\rho_k\sum_{j\in
 J_k}\lambda_{k,j}T(e_j\otimes\xi)\Vert^p\;\;\Sim^{(1+\eps)}\;\;
 \sum_{k=1}^n\vert\rho_k\vert^p\Vert\sum_{j\in
 J_k}\lambda_{k,j}T(e_j\otimes\xi)\Vert^p$$
where as usual the abbreviation $a\Sim\limits^C b$ means
$\max(a/b,b/a)\leq \sqrt C$. Let $S^{(n)}_\xi$ be the isometry
$\ell^p_n\hookrightarrow\ell^p$ defined by
 $$S^{(n)}_\xi (e_k)=\sum_{j\in J_k}\lambda_{k,j} e_j$$
and let
 $$T_\xi^{(N,n)}=T_N\rond(S^{(n)}_\xi\otimes Id_F):\ell_p^n(F)\to X$$
Then clearly for $N$ sufficiently large $T^{(N,n)}_\xi$ is a
$K$-isomorphic embedding and
 $$\lim_{N,\U}\Vert T^{(N,n)}_\xi
 (\sum_{k=1}^n\rho_ke_k\otimes \xi)\Vert_X^p\;\;\Sim^{1+\eps}\;\;
 \lim_{N,\U}\sum_{k=1}^n\vert\rho_k\vert^p
 \Vert T^{(N,n)}_\xi(e_k\otimes \xi)\Vert_X^p$$
for all $\rho_1,...,\rho_n\in \comp$. Using the compacity of the
unit ball of $\ell_p^n$ one easily deduces that for some
$N=N(n,\eps)$
 $$\Vert T^{(N,n)}_\xi(\sum_{k=1}^n\rho_ke_k\otimes \xi)\Vert_X^p
 \;\;\Sim^{1+2\eps}\;\; \sum_{k=1}^n\vert\rho_k\vert^p
 \Vert T^{(N,n)}_\xi(e_k\otimes \xi)\Vert_X^p$$
for all $\rho_1,...,\rho_n\in \comp$. Set
$T_\xi^{(n)}=T_\xi^{(N(n,\eps),n)}$. Note that
 $$\Vert T_\xi^{(n)}(\sum_{k}u_k\otimes\xi)\Vert_{X}^p
 \;\;\Sim^{(1+2\eps)^2}\;\;
 \sum_{k}\Vert T_\xi^{(n)}(u_k\otimes \xi)\Vert_{X}^p$$
for every sequence $(u_k)$ of pairwise disjoint elements of 
$\ell_p^n$.

\noindent Note moreover that if some $x\in F$ verifies
 $$\Vert T_n(\sum_{k}u_k\otimes x)\Vert_{X}^p
 \;\;\Sim^{1+\delta}\;\;
 \sum_{k}\Vert T_n(u_k\otimes x)\Vert_{X}^p$$
for every $n$ and every disjoint sequence $(u_k)$ of $\ell_p^n$,
then we have for every $n$:
 $$\Vert T_\xi^{(n)}(\sum_{k}u_k\otimes x)\Vert_{X}^p
 \;\;\Sim^{1+\delta}\;\;
 \sum_{k}\Vert T_\xi^{(n)}(u_k\otimes x)\Vert_{X}^p$$
too (simply because
the $S_\xi^N(u_k)$ are pairwise disjoint too). Let us say that the
sequence $(T_\xi^{(n)})$ is a $(\eps,\xi)$-refinement of the
sequence $(T_n)$.

Let now   $\E=\{\xi_1, ..., \xi_d\}$ be an $\eps$-net in the unit
sphere $S_F$ of $F$, i.e. a maximal set of $\eps$-separated points
of $S_F$. We define successively sequences $(T_{\xi_1}^{(n)}),
(T_{\xi_1,\xi_2}^{(n)})$, ..., $(T_{\xi_1, ..., \xi_d}^{(n)})$ of
$K$-isomorphic embeddings of $\ell_p^n(F)$ into $X$. The sequence
$((T_{\xi_1}^{(n)}))$ is a $(\eps,\xi_1)$-refinement of the
sequence $(T_n)$, and for every $j=2,...,d$, the sequence
$(T_{\xi_1,...,\xi_j}^{(n)})$ is a $(\eps,\xi_j)$-refinement of
the sequence $(T_{\xi_1,...,\xi_{(j-1)}}^{(n)})$. The final
operators, which we  denote by $T_\E^{(n)}$, are still
$K$-isomorphic embeddings and  verify
 $$\forall\xi\in \E,\forall(\rho_k)\in\ell_p^n,\quad
 \Vert T_\E^{(n)}(\sum_k\rho_ke_k\otimes \xi)\Vert_X^p
 \;\;\Sim^{(1+2\eps)^2}\;\;
 \sum_k\vert\rho_k\vert^p\Vert T_\E^{(n)}
 (e_k\otimes \xi)\Vert_X^p$$
If now $x\in S_F$ is arbitrary, let $\xi\in \E$ with
$\N{x-\xi}\leq\eps$. For every norm one $(\rho_k)\in\ell_p$ we
have by triangular inequality in $X$ (in the Banach case):
 $$\left\vert\,\Vert\sum_k\rho_kT_\E^{(n)}
 (e_k\otimes x)\Vert-\Vert\sum_k\rho_kT_\E^{(n)}(e_k\otimes
 \xi)\Vert\,\right\vert
 \leq\Vert T_\E^{(n)}\Vert\,\Vert
 \sum_k\rho_k e_k\otimes (x-\xi)\Vert_{\ell_p(F)}
 \leq \eps K$$
and similarly by triangular inequality in $\ell_p^n$, and in $X$:
 $$\left\vert\,(\sum_k\vert\rho_k\vert^p
 \Vert T_\E^{(n)}(e_k\otimes x)\Vert^p)^{1/p}
 -(\sum_k\vert\rho_k\vert^p\Vert T_\E^{(n)}
 (e_k\otimes \xi)\Vert^p)^{1/p}\right\vert\leq \eps K$$
Then we deduce that $(T_\E^{(n)}(e_k\otimes x))_k$ is $f(\eps,
K)$-equivalent to the $\ell_p^n$-basis, with $f(\eps, K)\to 1$
when $\eps\to 0$  (we find $\ds f(\eps, K)\leq
{(1+2\eps)^4(1+2K\eps)\over 1-2K\eps(1+2\eps)^2}$). Similar
estimations hold in the $p$-normed case.

A careful examination of what has been done shows that each
$T_\E^{(n)}$ is deduced from some $T_N$ simply by composing on the
right with some $S\otimes {\rm Id}_F$, where $S: \ell_p^n\hookrightarrow
\ell_p^N$ is an isometry. Note that $S$ maps the basis vectors of
$\ell_p^n$ onto disjoint vectors in $\ell_p^N$ (if $p=2$ it is not
automatic but results from the construction). If
$p\geq 1$ the range of such an isometry is always
$1$-complemented by some norm one projection $Q_S$ (see [LT],
prop. 2.a.1; in fact this projection verifies $\N{Q_S(\abs{x})}_p\leq
\N{x}_p$ for every $x\in\ell_p^N$). Then
$\widetilde Q_S=Q_S\otimes {\rm Id}_F$ is a norm
one projection in $\ell_p^N(F)$. If $P_N$ is a projection from $X$ onto the
range of $T_N$, then $T_N\widetilde
Q_ST_N^{-1}P_N$ is a projection  from $X$ onto the range of $T_\E^{(n)}$
(with norm $\leq K\N{P_N}$).\cqfd

\smallskip

\proclaim Lemma 3.3. For every $j\geq 1$ let $\U_j$ be a free
ultrafilter over $I$ and $(X_{i,j})_{i\in I}$ be a family of
$d_j$-dimensional subspaces of $L_p(\rA)$ such that
$(\prod_iX_{i,j})_{\U_j}$, considered as a subspace of
$L_p(\rA)_{\U_j}$,  is disjoint from $L_p(\rA)$. Let $(\eps_j)$ be
an arbitrary sequence of positive real numbers. Then there exist a
sequence $(i_j)_j$ in I and a sequence $(p_j)$ of pairwise
disjoint projections of $\rA$ such that:\hfill\break
\centerline{$\forall j\geq 1,\quad\sup\{\N{h-p_jhp_j}\mid h\in
X_{i_j,j};\N{h}\leq 1\}\leq \eps_j$.}

\proof Given $j$, a finite system $\P$ of pairwise disjoint
projections, and a finite dimensional subspace $V$ of $L_p(\rA)$,
we can obtain, using Lemma 2.6,
 a family $(s_i)_{i\in I}$ of projections of $\rA$ which
commute with $\P$ and such that:\medskip

i) $\forall k\in V,\qquad\ds\N{s_ik}+\N{ks_i}\To_{i,\U_j}0$

ii) for every bounded family $(h_i)\in \prod_i X_{i,j}$,
$\N{s_i^\perp h_i}+\N{h_is_i^\perp }\To\limits_{i,\U_j}0$.

\noindent To see this we just note that $V$ and
$(\prod_iX_{i,j})_{\U_j}$  are separable since they are finite
dimensional.

By compacity of the unit balls of $V$ and of $\prod_i X_{i,j}$,
the conditions (i), (ii) clearly imply:\medskip

i') $\ds\sup\{\N{s_ik}+\N{ks_i}\mid k\in V,\N{k}\leq
1\}\To_{i,\U_j}0$.

ii') $\sup\{\N{s_i^\perp h}+\N{hs_i^\perp}\mid h\in X_{i,j},
\N{h}\leq 1\}\To\limits _{i,\U_j}0$.

Let $(\delta_j)$ be a sequence of positive real numbers. Now we
construct by induction a sequence $(i_j)$ of distinct indices in
$I$ and a sequence of pairwise commuting projections $(q_j)$ such
that for every $j\ge 1$
 $$\eqalign{\forall h\in \sum_{n\leq j-1}X_{i_n,n},
 \quad
 &\N{q_jh}+\N{hq_j}<\delta_j\N{h}\cr
 \forall h\in X_{i_j,j},\quad
 &\N{q_j^\perp
 h}+\N{hq_j^\perp}<\delta_j\N{h}}$$
We shall consider the convex case ($p\geq 1$), the $p$-normed case
($0<p<1$) being treated analogously. Choose some  $i_1\in I$ and
set $q_1=\1$.  Assume constructed $i_1,...,i_j$ and $q_1,...,q_j$.
Set $V=\sum\limits_{n=1}^jX_{i_n,n}$, and let  $(s_i)$ be a family
verifying the conditions (i'), (ii') above with $j+1$ in place of
$j$. Thus for some $i\in T\setminus\{i_1,...,i_j\}$ we have:
 $$\eqalign{\forall h\in \sum_{n\leq j}X_{i_n,n},\quad
 &\N{s_ih}+\N{hs_i}<
  \delta_{j+1}\N{h}\cr
 \forall h\in X_{i,j+1},\quad
 &\N{s_i^\perp
 h}+\N{hs_i^\perp}<\delta_{j+1}\N{h}}$$
Then set $i_{j+1}=i$ and $q_{j+1}= s_i$. Finally, define
$p_j=q_j\bigwedge\limits_{k>j} q_k^\perp$. It is easy to check
that the two sequences $(i_j)$ and $(p_j)$ satisfy the
requirements of the lemma if the $\delta_j$ are sufficiently
small. \cqfd

\smallskip

 Now we are in a position to prove Theorem 3.1.\medskip
\smallskip

\noindent {\bf Proof of Theorem 3.1:} By Lemma 3.2, we find
$K$-embeddings $T_{j,n}$ of $\ell_p^n(F_j)$ into $X$ such that for
each nonzero $x\in F_j$ the sequence
$\ds\big(\N{T_{j,n}(e_i\otimes x)}^{-1}T_{j,n}(e_i\otimes
x)\big)_{1\leq i\leq n}$ is $(1+1/n)$-equivalent to the
$\ell_p^n$-basis. Let $\U$ be a free ultrafilter over $\nat$ and
define $\widetilde
T_j:\bigcup\limits_n\ell_p^n(F_j)\hookrightarrow X_\U$ by
$\widetilde T_j(e_i\otimes x)=(T_{j,n}(e_i\otimes x))^\bullet$,
with the agreement that $T_{j,n}(e_i\otimes x)=0$ if $i> n$. Then
$\widetilde T_j$ is a $K$-embedding into $L_p(\A)=L_p(\rA)_\U$,
that we extend by continuity to the whole of $\ell_p(F_j)$. For
every nonzero $x\in F_j$, the sequence $\ds\big(\Vert\widetilde
T_j(e_i\otimes x)\Vert^{-1}\widetilde T_j(e_i\otimes x)\big)_i$ is
1-equivalent to the $\ell_p$-basis, so it defines in $L_p(\A)_\U$
an element  disjoint from $L_p(\A)$, and a fortiori from
$L_p(\rA)$. We can identify $L_p(\A)_\U=\big(L_p(\rA_\U)\big)_\U$
with $L_p(\rA)_{\U\times\U}$. Set $S_{j,i,n}: F_j\to X, x\mapsto
T_{j,n}(e_{i}\otimes x)$; for $n\geq i$ the operator $S_{j,i,n}$
induces a $K$-isomorphic embedding of $F_j$ into $X$. Moreover, if
the initial copies of $\ell_p^n(F_j)$ in $X$ are $C$-complemented
in $L_p(\rA)$,  the ranges $S_{j,i,n}(F_j)$, $n\geq i$, are
$KC$-complemented. For every $x\in F_j$, the double sequence
$(S_{j,i,n}(x))_{i,n}$ defines an element of
$L_p(\rA)_{\U\times\U}$ disjoint from $L_p(\rA)$.

Let $\V$ be the trace of the ultrafilter $\U\times \U$ over the
set $D=\{(i,n)\mid n\geq i\}$. We apply Lemma 3.3 to the family
$(S_{j,i,n}(F_j))_{(i,n)\in D}$ and the ultrafilter $\V$.  Let
$(i_j,n_j)$ be a sequence in $D$ and $(p_j)$  be a disjoint
sequence of projections of $\rA$ satisfying the conclusion of that
lemma. From now on we write $E_j$ in place of
$S_{j,i_j,n_j}(F_j)$: recall that $E_j$ is $K$-isomorphic to $F_j$
by some isomorphism $T_j: F_j\to E_j$; and that if we denote by $R_j$
the operator $L_p(\rA)\to L_p(\rA)$, $h\mapsto p_jhp_j$ we have
$\Vert (Id-R_j)\mid_{E_j}\Vert<\eps_{j}$, which proves assertion
ii) of the theorem. In the case iii), we have moreover that $E_j$
is $KC$-complemented in $L_p(\rA)$ by some projection $P_j$.

Now we can easily accomplish the proof of the theorem. The
assertion i) of the theorem follows by a standard perturbation
argument which we sketch here (in the convex case) for further use
in Section 6. Let $(p_j)$ and $(E_j)$ be as before. Then for every
finite sequence $(y_j)\in\prod_j E_j$ we have
 $$\Vert\sum_jy_j-p_jy_jp_j\Vert
 \leq\sum_j\eps_j\Vert y_j\Vert
 \leq\sum_j{\eps_j\over 1-\eps_j} \Vert p_jy_jp_j\Vert
 \leq\eps\sup_j\Vert p_jy_jp_j\Vert$$
where $\ds\eps=\sum_{j\geq 1}{\eps_j\over 1-\eps_j}$ is finite and
small if the $\eps_j$'s are sufficiently small. On the other hand,
since the projections $p_j$ are pairwise disjoint,
 $$\Vert\sum_j p_jy_jp_j\Vert=(\sum_j\N{p_jy_jp_j}^p)^{1/p}$$
Thus it follows that:
 $$(1-\eps)(\sum_j\N{p_jy_jp_j}^p)^{1/p}
 \leq\Vert \sum_j y_j\Vert
 \leq(1+\eps)(\sum_j\N{p_jy_jp_j}^p)^{1/p}$$
However
 $$\N{p_jy_jp_j}\leq\N{y_j}\quad \hbox{ and}\quad
 \N{p_jy_jp_j}\geq(1-\eps_j)\N{y_j}\geq (1-\eps)\N{y_j}$$
Hence:
 $$(1-\eps)^2(\sum_j\N{y_j}^p)^{1/p}
 \leq\Vert \sum_j y_j\Vert\leq
 (1+\eps)(\sum_j\N{y_j}^p)^{1/p}$$
Assume now w.l.o.g. that $\N{T_j^{-1}}\leq 1$, $\N{T_j}\leq K$ for
every $j\geq 1$. We define
 $$T: F=(\mathop\bigoplus\limits_{j\geq 1} F_j)_p\to E=
 \sum\limits_{j\geq 1}E_j\;\hbox{ by }\;
 T((x_j))=\sum\limits_j T_j x_j$$
Then from the preceding inequalities we deduce that
$\N{T^{-1}}\leq (1-\eps)^{-2}$, $\N{T}\leq (1+\eps)K$. This proves
assertion i).

In case iii) of the theorem, since
 $$\Vert(Id-P_jR_j)\mid_{E_j}\Vert
 =\Vert(P_j-P_jR_j)\mid_{E_j}\Vert<KC\eps_j$$
it follows that  $W_j=P_jR_j\mid_{E_j}$ is for small $\eps_j$ an
isomorphism $E_j\to E_j$, with $\N{W_j^{-1}}\leq
(1-KC\eps_j)^{-1}$. Then   $Qx =\sum_j W_j^{-1}P_jR_jx$ defines a
bounded projection from $X$ onto $\sum_j E_j$. In fact
$\N{Qx}_p^p\leq
[(1+\eps)(1-\eps)^{-2}]^p\sum\limits_j\N{W_j^{-1}P_jR_jx}^p\leq
M^p \sum\limits_j\N{R_jx}^p=M^p\N{x}^p$ with
$M=KC(1+\eps)(1-\eps)^{-2}(1-KC\eps)^{-1}$. \cqfd

\bg\centerline{\bf 4. Equiintegrability and the Subsequence
Splitting Lemma}
\medskip

 In [R3], N. Randrianantoanina introduced the notion of
$p$-equiintegrable susbset of a non-commutative $L_p$-space. We
give here a seemingly more restrictive definition of
$p$-equiintegrable sets; it will appear later that this definition
is in fact equivalent to Randrianantoanina's one.

\proclaim Definition {4.1}. Let $\rA$ be a  von Neumann algebra
and $0<p<\infty$. A bounded subset $K$ of $L_p(\rA)$ is called
{\rm $p$-equiintegrable} if $\sup\limits_{h\in K}\N{e_\alpha
he_\alpha}_p\To\limits_\alpha 0$ for every net $(e_\alpha)$ of
projections of $\rA$ which w*-converges to $0$.

\noindent{\bf Remark 4.2.} {\it Finite subsets of $L_p(\rA)$ are
$p$-equiintegrable.} In fact, given a net of projections
$(s_\alpha)$ w*-converging to $0$,  let $A$ be the set of positive
reals $p$ such that $\Vert hs_\alpha\Vert_p\To\limits_\alpha 0$
for every $h\in L_p(\rA)$.  By H\"{o}lder's inequality, one easily
sees that $q\in A$ whenever $0<q<p$ and $p\in A$. Thus  $A$ is an
interval whose left endpoint is 0.  On the other hand, if $p\in
A$, then $2p\in A$ for
 $$\Vert hs_\alpha\Vert_{2p}^2=\Vert s_\alpha h^*hs_\alpha\Vert_p
 \leq \Vert(h^*h)s_\alpha\Vert_p$$
However, it is clear that  $2\in A$ since $\N{h s_\alpha}_2^2=\langle
h^*h,s_\alpha\rangle$ in the identification of $L_1(\rA)$ with $\rA_*$.
Therefore, we deduce that
$A=(0,\infty)$.

\proclaim Lemma {4.3}. Every net $(p_\alpha)$ of $\sigma$-finite
projections of $\rA$ which w*-converges to $0$ contains a sequence
which still w*-converges to $0$.

\proof Construct inductively a sequence $(\ph_n)$ in $\rA_*^+$ and
a sequence $(\alpha_n)$ such that:

i) $\max\,\{\ph_m(p_{\alpha_n})\mid m=1,...,n-1\}<n^{-1}$

ii) $\ph_n$ has support $p_{\alpha_n}$ and norm 1.

\noindent Set $\psi=\sum_{m=1}^\infty 2^{-m}\ph_m$, then $s(\psi)$
dominates all the $p_{\alpha_n}$'s, and clearly
$\psi(p_{\alpha_n})\To\limits_{n\to\infty}0$. Hence
$(p_{\alpha_n})$ w*-converges to zero. \cqfd

\proclaim Proposition {4.4}. A bounded subset $K$ of $L_p(\rA)$ is
$p$-equiintegrable if and only if for every sequence $(e_n)$ of
projections of $\rA$ which w*-converges to $0$ we have
$\sup\limits_{h\in K}\N{e_n he_n}_p\To\limits_n 0$. In particular,
a subset $K$ of $L_p(\rA)$ is $p$-equiintegrable if every
countable subset of $K$ is.

\proof The condition is clearly necessary. Conversely if $K$ is
not equiintegrable, there exists a family $(p_i)_{i\in I}$ of
projections of $\rA$ and an ultrafilter $\U$ over $I$ such that
w*-$\lim\limits_{i,\U} p_i=0$ but
$\lim\limits_{i,\U}\sup\limits_{h\in K}\N{p_ihp_i}_p>\delta>0$. We
can clearly suppose that for some family $(h_i)$ of elements of
$K$ we have $\N{p_ih_ip_i}\geq\delta/2$ for every $i\in I$. Let
$p'_i=\ell(p_ih_ip_i)\vee r(p_i h_ip_i)$: then $p'_i\leq p_i$ and
$(p'_i)$ w*-converges to $0$ with respect to $\U$, each $p'_i$ is
$\sigma$-finite and $p'_ih_ip'_i=p_ih_ip_i$.  By Lemma 4.3, there
exists a subsequence $(p'_{i_n})$ which w*-converges to $0$.\cqfd

\smallskip
\noindent{\bf Remark 4.5.} {\it If $K$ is $p$-equiintegrable, then
 for all bounded nets $(x_\alpha)$, $(y_\alpha)$ of  positive
elements of $\rA$ which w*-converge  to 0, we have
$\sup\limits_{h\in K}\N{x_\alpha hy_\alpha}_p\To\limits_\alpha 0$.}

\proof  Fix $\eps>0$ and let $e_{\alpha,\eps}$ be the spectral
projection $\chi_{[\eps,+\infty)}(x_\alpha+y_\alpha)$. Since
$e_{\alpha,\eps}\leq \eps^{-1}(x_\alpha+y_\alpha)$, we have
$e_{\alpha,\eps}\To\limits^{w^*}_\alpha 0$, hence
$\sup\limits_{h\in K}\N{e_{\alpha,\eps}
he_{\alpha,\eps}}_p\To\limits_\alpha 0$. Consequently we have
$\sup\limits_{h\in K}\N{x_\alpha e_{\alpha,\eps}
he_{\alpha,\eps}y_\alpha}_p\To\limits_\alpha 0$. On the other
hand, since $0\leq x_\alpha\leq (x_\alpha+y_\alpha)$, there exist
$c_\alpha\in\rA$, $0\leq c_\alpha\leq \1$, such that
$x_\alpha=(x_\alpha+y_\alpha)^{1/2}
c_\alpha(x_\alpha+y_\alpha)^{1/2}$. Then
 $$\N{x_\alpha e_{\alpha,\eps}^\perp}
 \leq\N{(x_\alpha+y_\alpha)^{1/2}c_\alpha}\N{(x_\alpha
 + y_\alpha)^{1/2}e_{\alpha,\eps}^\perp}
 \leq \eps^{1/2}M^{1/2}$$
where $M$ is a bound for the $\N{x_\alpha+y_\alpha}$. Similarly,
$\N{e_{\alpha,\eps}^\perp y_\alpha}\leq \eps^{1/2}M^{1/2}$. So we
obtain $\adh{\lim\limits_\alpha}\sup\limits_{h\in K}\N{x_\alpha
hy_\alpha}_p\leq 2\eps^{1/2} M^{3/2}M'$, where $M'$ is a bound for
the $\N{h}$, $h\in K$.\cqfd

\smallskip

Now we  characterize the $p$-equiintegrability of a bounded
sequence in $L_p(\rA)$ in terms  of the element it defines in an
ultrapower of $L_p(\rA)$ and the disjointness of this element from
$L_p(\rA)$. To this end  we introduce the following notation.  Let
$\U$ be an ultrafilter over the index set $I$. Let $s_e$ be the
support of $L_p(\rA)$ in $L_p(\rA)_\U$ (considered as a
non-commutative $L_p$-space $L_p(\A)$). We have thus (since
$L_p(\rA)$ is self-adjoint and generated by its positive cone):
 $$\eqalign{s_e
 &=\sup\{\ell(\hat h)\vee r(\hat h)\mid h\in L_p(\rA)\}\cr
 &=\sup\{\ell(\hat h)\mid h\in L_p(\rA)\}
 =\sup\{r(\hat h)\mid h\in L_p(\rA)\}\cr
 &=\sup\{s(\hat h)\mid h\in L_p(\rA)_+\}}$$
It is also clear that $s_e$ does not depend on $p\in(0,\infty)$,
since $s(\hat h)=s((\hat h)^p)=s(\widehat{h^p})$ for every $h\in
L_p(\rA)^+$. Note also that an element $\t h\in L_p(\rA)_\U$ is
disjoint from $L_p(\rA)$ iff $\t h=s_e^\perp \t h s_e^\perp$.

If $\rA$ is $\sigma$-finite, then $s_e=s(\hat h_0)$ for every
$h_0\in L_p(\rA)_+$ with support $s(h_0)=\1$ (when $p=1$ this
means that the associated $\ph_0\in\rA_*$ is faithful). For, let
$h\in L_p(\rA)_+$; since $\rA\cdot h_0$ is dense in $L_p(\rA)$,
there exists for every $\eps>0$ an $x\in \rA$ such that $\N{h-x
h_0}<\eps$. Then $\Vert\hat h-\hat x\hat h_0\Vert\leq\eps$ and we
see that $\hat h$ is in the closure of $\A\hat h_0$. So $s(\hat
h)=r(\hat h)\leq r(\hat h_0)=s(\hat h_0)$.

If $\rA$ is a finite von Neumann algebra, then $s_e$ is a central
projection. For, assume that there is a finite normal faithful
trace  $\tau$  on $\rA$. Let $\hat\tau=(\tau)^\bullet$ be its
canonical image in $(\rA_*)_\U$: then $s_e=s(\hat \tau)$. For any
$\t x,\t y\in \rA_\U$  we clearly have $\hat\tau(\t x\t
y)=\hat\tau(\t y\t x)$. By the w*-density of $\rA_\U$ in $\A$, we
deduce that $\hat\tau$ is tracial, and consequently its support is
central. In the general case, we can argue similarly, using a
faithful family of normal traces with pairwise disjoint supports.

\smallskip
 The main result of this section is the following.
Recall that an ultrafilter $\U$ is {\it countably incomplete} if
there exists a sequence $(A_n)_{n\geq 1}$ of members of $\U$ such
that $\bigcap\limits_{n\geq 1}A_n=\emptyset$ (so is every non
trivial ultrafilter on a countable set).

\proclaim Theorem {4.6}. Let $\U$ be a countably incomplete
ultrafilter over the set $I$. Let $s_e$ be the support of
$L_p(\rA)$ in $L_p(\A)=L_p(\rA)_\U$. Then an element $\t h$ of
$L_p(\rA)_\U$ verifies the equality $s_e^\perp \t h s_e^\perp=0$
if and only if it admits a $p$-equiintegrable representing family
$(h_i)_{i\in I}$.

For the proof of this theorem we shall need the following density
lemma (see also [Ju2] and [JX] for similar results).

\proclaim Lemma 4.7. Let $h_0$ be an element of $L_p(\rA)$,
$0<p<\infty$. Then $\rA\cdot h_0$ is dense in $L_p(\rA)\cdot
r(h_0)$, and $h_0\cdot \rA$ is dense in $\ell(h_0)\cdot L_p(\rA)$.

\proof Clearly we can assume w.l.o.g. that $h_0$ is positive and
$\N{h_0}_p=1$. Then $\ell(h_0)=r(h_0)=s(h_0)$. Assume first that
$p\geq 1$ and let $q$ be the conjugate index of $p$. Then the dual
space of $L_p(\rA)\cdot s(h_0)$ is the space $s(h_0)\cdot
L_q(\rA)$ (under the duality $\langle h,k\rangle={\rm Tr\,}(hk)$).
If $k\in s(h_0)\cdot L_q(\rA)$ belongs to the annihilator of
$\rA\cdot h_0$, then ${\rm Tr}(xh_0k)=0$ for every $x\in\rA$,
which in turn implies that $h_0\cdot k=0$ (as element of
$L_1(\rA)$), hence $k=s(h_0)\cdot k=0$. So the linear space
$\rA\cdot h_0$ is dense in $L_p(\rA)\cdot s(h_0)$. Similarly,
$h_0\cdot \rA$ is dense in $s(h_0)\cdot L_p(\rA)$. Assume now that
$1/2\leq p<1$. Every $h\in L_p(\rA)\cdot s(h_0)$ can be factorized
as:
 $$h=u\abs{x}=u\abs{x}^{1/2}\abs{x}^{1/2}
 =(u\abs{x}^{1/2}s(h_0))\cdot (\abs{x}^{1/2}s(h_0))$$
since the supports of $\abs{h}$ and $\abs{h}^{1/2}$ coincide with
the right support of $h$, hence are included in $s(h_0)$. So
$h=k'k''$ with $k', k''\in L_{2p}(\rA)\cdot s(h_0)$. Since $2p\geq
1$, there exists by the preceding argument a sequence $(y_n)$ in $\rA$
such that $y_n\cdot h_0^{1/2}\to k''$ (for the norm of $L_{2p}(\rA)$);
and for every $n$ there exists a sequence $(x^n_m)_m$ in $\rA$
such that $x^n_m\cdot h_0^{1/2}\to k'y_ns(h_0)$ when $m\to\infty$.
Then
 $$\lim\limits_{n\to\infty}\lim\limits_{m\to\infty}x^n_mh_0
 =\lim\limits_{n\to\infty}k'y_ns(h_0)\cdot h_0^{1/2}=k'k''=h$$
Therefore the conclusion of the lemma is true for $1/2\leq p<1$.
Iterating this procedure we see that it is true for every
$0<p<1$.\cqfd

\smallskip
\noindent{\bf Proof of Theorem 4.6:} Let $(h_i)$ is a bounded
family in $L_p(\rA)$, and let $\t h$ be the corresponding element
of $L_p(\rA)_\U$.

\noindent a) Suppose that $s_e^\perp\t h s_e^\perp\neq 0$. Then
there exists   a $\sigma$-finite projection $q\in\A$ disjoint from
$s_e$ such that $q\t h q\neq 0$. Let $\Vert q\t h
q\Vert=\delta>0$. By Theorem 2.3, for every $\ph\in\rA_*^+$ there
exists a family $(q_i)$ of projections of $\rA $ such that $\t
q:=(q_i)^\bullet\geq q$ and $\t q^\perp\geq s(\widehat \ph)$. The
second inequality yields $\lim\limits_{i,\U}\ph(q_i)=0$, and the
first one implies $\lim\limits_{i,\U}\N{q_ih_iq_i}\geq \Vert q\t h
q\Vert=\delta$. Note that we may suppose that each $q_i$ is
$\sigma$-finite, by replacing if necessary $q_i$ by
$q'_i=\ell(q_ih_iq_i)\vee r(q_i h_i q_i)$. So for each $\eps>0$ we
can find $i\in I$ such that $\N{q_ih_iq_i}\geq \delta/2$ and
$\ph(q_i)<\eps$. Now it is easy to construct inductively a
sequence $(\ph_n)$ in $\rA_*^+$, a sequence $(i_n)$ in $I$, and a
sequence $(s_n)$ of $\sigma$-finite projections of $\rA$ such
that:\medskip

i) $\N{s_nh_{i_n}s_n}\geq \delta/2$

ii) $\max\,\{\ph_m(s_n)\mid m=1,...,n-1\}<\ds 1/n$

iii) $\ph_n$ has support $s_n$ and norm 1.\medskip

\noindent Set $\psi=\sum_{m=1}^\infty 2^{-m}\ph_m$, then $s(\psi)$
dominates all the $s_n$'s, and clearly
$\psi(s_n)\To\limits_{n\to\infty}0$. Hence $(s_n)$ w*-converges to
zero, and by (i) the sequence $(h_{i_n})$ cannot be
$p$-equiintegrable, and a fortiori the family $(h_i)$ is not
$p$-equiintegrable.

\noindent b) Conversely, assume that $s_e^\perp\t h s_e^\perp=0$.
Then $\t h=\t h s_e+ s_e(\t hs_e^\perp)\in L_p(\A)s_e+s_eL_p(\A)$.
By Lemma 4.7, $\A\cdot L_p(\rA)=\span\{x\hat k\mid x\in\A, k\in
L_p(\rA)\}$ is dense in $L_p(\A)\cdot s_e$, and $L_p(\rA)\cdot\A$
 dense in $s_e\cdot L_p(\A)$. Note also that by Proposition 2.1,
$\A\cdot L_p(\rA)=\rA_\U\cdot L_p(\rA)$ and $L_p(\rA)\cdot\A =
L_p(\rA)\cdot\rA_\U$.  If $\t h\in \rA_\U\cdot L_p(\rA)$, it
admits a representing family $(h_i)$ of the type $(x_ih)$, where
$(x_i)$ is a bounded family in $\rA$ and $h\in L_p(\rA)$. Let
$(s_\alpha)$ be a net of projections which w*-converges to $0$.
Then by  Remark 4.2,
 $$\sup_i\N{h_is_\alpha}_p=\sup_i\Vert x_ihs_\alpha\Vert_p
 \leq (\sup\limits_i\N{x_i})\,\Vert he_\alpha\Vert_p
 \To\limits_\alpha 0$$
Thus $(h_i)$ is $p$-equiintegrable. Similarly, every $\t h\in
L_p(\rA)\cdot \rA_\U$ has a $p$-equiintegrable representing
family. Hence the proof will be complete if we show that the
subspace of $L_p(\rA)_\U$ consisting of elements having a
$p$-equiintegrable representing family is closed.

Let $(\t h^{(n)})_n$ be a sequence in $L_p(\rA)_\U$ which
converges to an element $\t h$ and suppose that each $\t h^{(n)}$
admits a $p$-equiintegrable representing family $( h^{(n)}_i)_i$.
We may suppose that $\Vert \t h^{(n)}-\t h\Vert<1/n$. Let $(h_i)$
be a representing family for $\t h$. Since $\U$ is countably
incomplete, we can find a decreasing sequence $U_1\supset
U_2\supset...\supset U_n\supset...$ of members of $\U$ such that
$\bigcap\limits_n U_n=\emptyset$ and
 $$i\in U_n\implies \Vert h^{(n)}_i-h_i\Vert\leq \frac 1n$$
Set $h'_i=h^{(n)}_i$ if $i\in U_n\setminus U_{n+1}$ and $h'_i=0$
if $i\not\in U_1$. Then $\Vert h'_i-h_i\Vert_p< n^{-1}$ for every
$i\in U_n$, which proves that $(h'_i)^\bullet=(h_i)^\bullet=\t h$.
Fix $n\geq 1$ and $i\in U_1$. Let $m\geq 1$ such that $i\in
U_m\setminus U_{m+1}$. If $m\geq n$ we have
 $$\Vert h'_i-h^{(n)}_i\Vert=\Vert h^{(m)}_i-h^{(n)}_i\Vert
 \leq \Vert h_i-h^{(m)}_i\Vert+\Vert h_i-h^{(n)}_i\Vert
 \leq \frac 1m+\frac 1n\leq \frac 2n$$
Consequently
 $$\forall i\in U_1,\ \inf_{1\leq j\leq n}\Vert h'_i-h^{(j)}_i\Vert
 \leq \frac 2n\ \To_{n\to\infty}0$$
Using the fact that a finite union of $p$-equiintegrable sets is
$p$-equiintegrable, we easily see that the family $(h'_i)_{i\in
U_1}$ is $p$-equiintegrable. Since $\{h'_i\mid i\in I\setminus
U_1\}=\{0\}$ is clearly $p$-equiintegrable too, we are done.\cqfd

\smallskip\noindent{\bf Remark.} In the proof of Theorem 4.6 the
hypothesis that the ultrafiter is countably incomplete was not
used for the sufficiency of the condition.

\smallskip
 Theorem 4.6 permits one to easily recover the following
 Subsequence Splitting
 Lemma obtained by N. Randrianantoanina [R3].

\proclaim Corollary {4.8}. Let  $0<p<\infty$ and $\rA$ be a von
Neumann algebra. Let $(h_n)$ be a bounded sequence in $L_p(\rA)$.
Then there exists an increasing sequence $(n_k)$ of integers and a
sequence $(p_k)$ of pairwise disjoint projections in $\rA$ such
that the sequence $(h_{n_k}-p_kh_{n_k}p_k)$ is $p$-equiintegrable.
As a consequence, we have the splitting $h_{n_k}=h'_k+h''_k$,
where $(h'_k)$ is $p$-equiintegrable and $(h''_k)$ is disjoint.

\proof  Let $\U$ be a free ultrafilter over $\nat$ and $\t
h=(h_n)^\bullet$. Let $\t h''=s_e^\perp\t h s_e^\perp$ and $\t
h'=\t h-\t h''$. Since $s_e^\perp\t h's_e^\perp=0$, by Theorem 4.6
the element $\t h'$ admits a representing sequence $(h'_n)$ which
is $p$-equiintegrable. Let $h''_n=h_n-h'_n$. Then $(h''_n)^\bullet
= \t h''$. Since $\t h''$ is disjoint from $L_p(\rA)$, by Theorem
2.7 there is an increasing sequence $(n_k)$ of integers and a
disjoint sequence $(p_k)$ of projections such that $\Vert
h''_{n_k}-p_kh''_{n_k}p_k\Vert \To 0$ when $k\to\infty$. Since
$(h'_n)$ is $p$-equiintegrable, $\Vert p_kh'_{n_k}p_k\Vert\to 0$,
and thus it follows that
 $$h_{n_k}-p_kh_{n_k}p_k=
 h'_{n_k}-p_kh'_{n_k}p_k+[h''_{n_k}-p_kh''_{n_k}p_k]$$
is a perturbation of a $p$-equiintegrable sequence by a norm
vanishing sequence, so is $p$-equiintegrable too. \cqfd

\proclaim Corollary {4.9}. The following conditions are equivalent
for a bounded subset  $K$ of $L_p(\rA)$:\hfill\break
 i) $K$ is $p$-equiintegrable;\hfill\break
 ii) for every disjoint sequence $(p_k)$ of
 projections of $\rA$,
 $\lim\limits_{k\to\infty}\sup\limits_{h\in K}\N{p_k
 h p_k}\to 0$;\hfill\break
 iii)  for every sequence $(p_k)$ of
projections of $\rA$ which decreases to $0$,
 $\lim\limits_{k\to\infty}\sup\limits_{h\in K}\N{p_k h p_k}\to 0$.
 \hfill\break
If in addition $\rA$ is $\sigma$-finite and $\ph_0$ is a normal
faithful state on $\rA$, i)-iii) are equivalent to the
following:\hfill\break
 iv) $\lim\limits_{\eps\to
 0}\sup\,\{\N{ehe}\mid h\in K,e\in \rA\;
 \hbox{projection such that}\;\ph_0(e)\leq\eps\,\}=0.$

\proof It is clear that  condition (i) implies (ii) and (iii). The
equivalence between (ii) and (iii) is easy (see also [R3]). To
prove that (ii) implies (i), suppose that $K$ is not
$p$-equiintegrable. Then by Proposition 4.4, it contains a
sequence $(h_n)$ which is not $p$-equiintegrable, and by Corollary
4.8, we can suppose that $h_n=h'_n+p_nh_np_n$, where $(h'_n)$ is
$p$-equiintegrable and $(p_n)$ is a disjoint sequence of
projections of $\rA$. Then $\N{p_nh_np_n}$ does not converges to
$0$, otherwise $(h_n)$ would be $p$-equiintegrable. The
equivalence of (iv) and (i) is due to the fact that a net
$(e_\alpha)$ of projections w*-converges to zero iff
$\ph_0(e_\alpha)\to 0$. This in turn follows from the density of
$\rA\cdot\ph_0$ in $\rA_*$ and the fact that
$x\cdot\ph_0(e_\alpha)=\ph_0(e_\alpha x)\leq
\ph_0(e_\alpha)^{1/2}\ph_0(x^*x)^{1/2}$ for every $x\in
\rA$.\cqfd\medskip

\noindent{\bf Remark.} Randrianantoanina [R3] took the part iii)
of Corollary 4.9 as the definition of $p$-equiintegrability.

\smallskip

Like for the  left resp. right disjoint sequences (see Remark
2.8), we shall also consider the corresponding left resp. right
$p$-equiintegrable sets.

\proclaim Definition {4.10}. We call a bounded set $K$ of
$L_p(\rA)$ {\rm left $p$-equiintegrable}, resp. {\rm right
$p$-equiintegrable} if for every net $(e_\alpha)$ of projections
w*-converging to $0$ we have $\sup\limits_{h\in K}\N{e_\alpha
h}\To\limits_\alpha 0$, resp. $\sup\limits_{h\in K} \N{he_\alpha}
\To\limits_\alpha 0$. We say that $K$ is {\rm
$p$-biequiintegrable} if it is both left  and right
$p$-equiintegrable.

In this definition we may again w.l.o.g. replace nets by
sequences. Note that if $K$ consists of positive elements, then
$K$ is $p$-equiintegrable iff it is $p$-biequiintegrable since
$\N{e h}=\N{eh^{1/2}h^{1/2}}\leq\N{ehe}\N{h}$ for every $h\in
L_p(\rA)^+$ and $e$ projection of $\rA$. Thus the four notions of
$p$-equiintegrability coincide on subsets of $L_p(\rA)^+$. Note
also that $K$ is left $p$-equiintegrable iff $K^*=\{h^*\mid h\in
K\}$ is right $p$-equiintegrable; if $K$ is left (resp. \nobreak
 right) $p$-equiintegrable and $B$ is a bounded subset of $\rA$,
then  the set  $\{k\cdot x\mid k\in K, x\in B\}$ (resp. $\{x\cdot
k\mid k\in K, x\in B\}$) is left  (resp. \nobreak  right)
$p$-equiintegrable too; in particular, $K$ is left (resp. right )
$p$-equiintegrable iff $\abs{K^*}:=\{\abs{h^*}\mid h\in K\}$
(resp. $\abs{K}:=\{\abs{h}\mid h\in K\}$) is $p$-equiintegrable.
Finally, $K$ is $p$-biequiintegrable iff both $\abs{K}$ and
$\abs{K^*}$ are.

\smallskip
Theorem 4.6 can be refined in the following way:

\proclaim Proposition {4.11}. Let $\U$ be a countably incomplete
ultrafilter over the set $I$. Let $s_e$ be the support of
$L_p(\rA)$ in $L_p(\A)=L_p(\rA)_\U$. Let $\t h\in L_p(\rA)_\U$.
Then:\hfill\break
 i) $\t h\in L_p(\A)s_e$ iff it admits a
right $p$-equiintegrable representing family.\hfill\break
 ii) $\t h\in s_eL_p(\A)$ iff it admits a
left $p$-equiintegrable representing family.\hfill\break
 iii) $\t h\in s_eL_p(\A)s_e$ iff it admits a
$p$-biequiintegrable representing family.

\proof If $\t h\not\in L_p(\A)s_e$,  $\t h s_e^\perp\neq 0$. Then
we can prove that $\t h$ has no right $p$-equiintegrable
representing family in a way very similar to the first part of the
proof of Theorem 4.6. For the converse implication we note that
$\rA_\U L_p(\rA)$ is dense in $L_p(\A)s_e$ (Lemma 4.7) and the
space of elements $\t h$ admitting a right $p$-equiintegrable
representing family is closed. The second assertion follows by
conjugation. Finally if $\t h$ admits a $p$-biequiintegrable
representing family, it belongs to $L_p(\A)s_e\cap
s_eL_p(\A)=s_eL_p(\A)s_e$. For the converse implication, we need
only to note that by Lemma 4.7 $L_p(\rA)\rA_\U L_p(\rA)$ is dense
in $s_eL_p(\A)s_e$ and that the space of elements $\t h$ admitting a
$p$-biequiintegrable representing family is closed.\cqfd

\smallskip
 The following corollary improves the Subsequence
Splitting Lemma.

\proclaim Corollary {4.12}. Let $0<p<\infty$ and $\rA$ be a von
Neumann algebra. Let $(h_n)$ be a bounded sequence in $L_p(\rA)$.
Then there exist an increasing sequence $(n_k)$ of integers and a
disjoint sequence $(p_k)$ of projections in $\rA$ such
that:\hfill\break
 - the sequence $(p_k^\perp h_{n_k}p_k^\perp)$ is
$p$-biequiintegrable,\hfill\break
 - the sequence $(p_k^\perp
h_{n_k}p_k)$ is left $p$-equiintegrable,\hfill\break
 - the sequence $(p_k h_{n_k}p_k^\perp)$ is
right $p$-equiintegrable.\hfill\break Consequently, the sequence
$(h_{n_k})$ splits into the sum of four sequences:\hfill\break
\centerline{$h_{n_k}=a_k+b_k+c_k+d_k$} where $(a_k)$ is
$p$-biequiintegrable, $(b_k)$ is left $p$-equiintegrable and
 right disjoint, $(c_k)$ is right $p$-equiintegrable and
left disjoint, $(d_k)$ is disjoint.

\proof Let $\U$ be a free ultrafilter over $\nat$. By Remark 2.8,
for every bounded sequence $(h_n)$ in $L_p(\rA)$ defining an
element  $\t h$ of $L_p(\A)s_e^\perp$ there exist a subsequence
$(h_{n_k})$ and a disjoint sequence of projections $(p_k)$ such
that $\N{h_{n_k}-h_{n_k}p_k}_p\To 0$. Moreover, given a finite set
of bounded sequences $(h_n^{(j)})$, $j=1,...,N$, each defining an
element of  $L_p(A)s_e^\perp$, one can find a {\it common}
increasing sequence $(n_k)$ and a {\it common} disjoint sequence
$(p_k)$ of projections such that $\Vert
h_{n_k}^{(j)}-h_{n_k}^{(j)}p_k\Vert_p\To 0$, $j=1,...,N$ (compare
with Lemma 3.3).

Now fix  a bounded sequence $(h_n)$ in $L_p(\rA)$, and le $\t h$
be the corresponding element in $L_p(\rA)_\U=L_p(\A)$. According
to the decomposition
 $$\t h=s_e\t hs_e + s_e\t hs_e^\perp
 + s_e^\perp\t hs_e + s_e^\perp\t hs_e^\perp$$
and by Proposition 4.11, we find four bounded sequences $(a_n)$,
$(d_n)$, $(c_n)$ and $(d_n)$ such that $h_n=a_n+b_n+c_n+d_n$,
$(a_n)$ $p$-biequiintegrable and  $(a_n)^\bullet=s_e\t h s_e$,
$(b_n)$ left $p$-equiintegrable and  $(b_n)^\bullet=s_e\t h
s_e^\perp$, $(c_n)$   right $p$-equiintegrable and
$(c_n)^\bullet=s_e^\perp\t h s_e$, and finally
$(d_n)^\bullet=s_e^\perp \t h s_e^\perp$.

Applying the preceding remark to the set
$\{(b_n),(c_n^*),(d_n),(d_n^*)\}$, we obtain an increasing
sequence $(n_k)$ of integers and a disjoint sequence $(p_k)$ of
projections such that the four sequences
 $$(b_{n_k}-b_{n_k}p_k),\ (c_{n_k}-p_kc_{n_k}),\
 (d_{n_k}-d_{n_k}p_k),\ (d_{n_k}-p_kd_{n_k})$$
all converge to $0$. Note that $(p_kb_{n_k})$ and $(c_{n_k}p_k)$
converge to zero too since $(b_n)$ and $(c_n)$ are respectively
left and right $p$-equiintegrable and the projections $p_k$ are
pairwise disjoint. Therefore, we deduce that the three sequences
$(b_{n_k}-p_k^\perp b_{n_k}p_k)$, $(c_{n_k}-p_kc_{n_k}p_k^\perp)$
and $(d_{n_k}-p_kd_{n_k}p_k)$ converge to zero as well. Thus we
can decompose $h_{n_k}$ as:
 $$h_{n_k}=a'_k+p_k^\perp
 b_{n_k}p_k+p_kc_{n_k}p_k^\perp+p_kd_kp_k$$
where $(a_k')$ is a  sequence such that
$\N{a'_k-a_{n_k}}\To\limits_{k\to\infty}0$. Consequently, $(a'_k)$
is $p$-biequiintegrable too. It follows that $(p_k^\perp
h_{n_k}p_k^\perp)=(p_k^\perp a'_kp_k^\perp)$ is
$p$-biequiintegrable, $(p_k^\perp h_{n_k}p_k)=(p_k^\perp
a'_kp_k+p_k^\perp b_{n_k}p_k)$ is left $p$-equiintegrable and
$(p_k h_{n_k}p_k^\perp)=(p_ka'_kp_k^\perp +p_kc_{n_k}p_k^\perp)$
is right $p$-equiintegrable.\cqfd

\smallskip \noindent{\bf Remark.} Using Corollary 4.12, one easily
sees that Corollary 4.9 extends to the left, right
$p$-equiintegrability and $p$-biequiintegrability.

\smallskip \noindent{\bf Remark.} If $\rA$ is finite, the notions of
left $p$-equiintegrability and right $p$-equiintegrability
coincide (so the four notions of $p$-equiintegrability coincide).
This can be deduced simply from Proposition 4.11 and the fact that
$s_e$ is central in this case. Consequently, in this case, the
sequences $(p_k^\perp h_{n_k}p_k)$ and $(p_k h_{n_k}p_k^\perp)$ in
Corollary 4.12 converge to zero.

\smallskip

The following gives one more characterization of
equiintegrability, which is quite useful in some context.

\proclaim Proposition {4.13}. Let  $0<p<\infty$ and $\rA$ be a von
Neumann algebra with unit ball $B_\srA$. Let $K$ be a subset of
$L_p(\rA)$. Then: \hfill\break
 i) $K$ is left  (resp. right)
$p$-equiintegrable iff for every $\eps>0$ there exists $h_\eps\in
L_p(\rA)$  such that for every $h\in K$, the distance
$d(h,h_\eps.B_\srA)$ (resp. $d(h,B_\srA\cdot h_\eps)$) in
$L_p(\rA)$ is majorized by $\eps$. \hfill\break
 ii) $K$ is $p$-equiintegrable iff  for every
$\eps>0$ there exists $h_\eps\in L_p(\rA)$  such that for every
$h\in K$, $d(h,B_\srA.h_\eps+h_\eps\cdot B_\srA)<\eps$.
\hfill\break
 iii) $K$ is $p$-biequiintegrable iff for every $\eps>0$ there
exists $h_\eps\in L_p(\rA)$  such that for every $h\in K$,
$d(h,B_\srA\cdot h_\eps\cdot B_\srA)<\eps$. \hfill\break
 In the case where
$\rA$ is $\sigma$-finite, one can take elements $h_\eps$ of the
form $M_\eps h_0$, where $M_\eps$ is a positive real number and
$h_0$ is a fixed positive element of $L_p(\rA)$ with full support
($s(h_0)=\1$).

\proof We give the proof for left equiintegrable sets, and a non
$\sigma$-finite von Neumann algebra; the other cases can be
treated similarly. The sufficiency of the condition is clear since
for every $h_0\in L_p(\rA)$, the set $h_0\cdot B_\srA$ is left
$p$-equiintegrable. Conversely, assume that $K$ is left
$p$-equiintegrable. Suppose that for some $\eps>0$ and for every
finite subset $F$ of $L_p(\rA)$ there exists $h_F\in K$ such that
$d(h_F,\sum\limits_{f\in F}(f\cdot B_\srA))>\eps$. Let $\cal F$ be
the net of finite subsets of $L_p(\rA)$, ordered by inclusion; let
$\Phi$ be the filter of cofinal subsets of $\cal F$ (generated by
the set of final sections $\Sigma_F=\{G\in {\cal F}\mid F\subset
G\}$); let finally $\U$ be an ultrafilter containing $\Phi$. By
Proposition 4.11, the element $\t h=(h_F)^\bullet$ of $L_p(\rA)_\U$
belongs to $s_e\cdot L_p(\A)$. Hence by Lemma 4.7, there exists
$h_\eps\in L_p(\rA)$ such that $d(\t h, \hat h_\eps\cdot
B_{\A})<\eps$. By Kaplansky's density theorem we have in fact
$d(\t h, \hat h_\eps\cdot B_{\srA_\U})<\eps$. Consequently, the set
$\{F\in{\cal F}\mid d(h_F, h_\eps\cdot B_\srA)<\eps\}$ belongs to
$\U$. Since the set $\{F\in {\cal F}\mid h_\eps\in
F\}=\Sigma_{\{h_\eps\}}$ belongs to $\U$ too, there is
$F\in\cal F$ such that $h_\eps\in F$ and $d(h_F, h_\eps\cdot
B_\srA)<\eps$, which contradicts the choice of the $h_F$'s. So in fact
for every
$\eps>0$ there exists $F_\eps=\{h_1^{(\eps)},...,h_n^{(\eps)}\}\in
\cal F$ such that $d(h,\sum\limits_{i=1}^n h_i^{(\eps)}\cdot
B_\srA)\leq\eps$ for every $h\in K$. Let $h_\eps =(\sum_{i=1}^n
h_i^{(\eps)}h_i^{(\eps)*})^{1/2}$, then for every $i=1,...,n$ we
have $h_i^{(\eps)}=h_\eps x_i$,  for some $x_i\in B_\srA$.
Consequently, $\sum\limits_{i=1}^n h_i^{(\eps)}\cdot B_\srA\subset
(nh_\eps)\cdot B_\srA$, and $d(h,(nh_\eps)\cdot B_\srA)\leq\eps$
for every $h\in K$.\cqfd\bigskip

\smallskip\noindent{\bf Historical comments.} i) In the
case of commutative $L_1$-spaces,  Corollary 4.8  was proved in
[KP] (where it is not explicitly stated but is a key ingredient of
the proof of the main result there). There are various extensions
to the Banach lattice setting; a general statement was given by L.
Weis using ultrapower techniques [W]. \hfill\break
 ii) In the non-commutative case a subsequence splitting lemma
 similar to Corollary 4.8 was obtained in [S] for symmetric spaces of
measurable
operators $E(\rA,\tau)$ associated with an order continuous
rearrangement invariant space $E$ and a von Neumann algebra $\rA$
equipped with a finite trace $\tau$ (see Lemma 1.1 and Proposition
2.2 of [S]).  Randrianantoanina proved Corollary 4.8 for symmetric
spaces of measurable operators $E(\rA,\tau)$ when $\tau$ is
semi-finite ([R1]) and for general non-commutative $L_p$-spaces
([R3]).\hfill\break
 iii) In the case of finite and $\sigma$-finite von Neumann
 algebras Proposition 4.13 goes back to [HRS].

\bigskip\noindent{\bf Application: weakly relatively compact
sets in $\rA_*$}
\medskip

  The $1$-equiintegrable
sets coincide with the weakly relatively compact sets. Proposition
4.13 can be used to give a new proof of some well known results of
C. A. Akeman ([A]):

\proclaim Theorem 4.14. Let $K$ be a bounded subset of the predual
$\rA_*$ of a von Neumann algebra $\rA$. The following assertions
are equivalent:\hfill\break
 i) $K$ is weakly relatively compact;\hfill\break
 ii) For every sequence $(p_n)$ of pairwise
disjoint projections in $\rA$,
$\lim\limits_{n\to\infty}\sup\limits_{\ph\in K}\vert\ph(p_n)\vert
=0$;\hfill\break
 iii) $K$ is 1-equiintegrable;\hfill\break
  iv) There exists $\psi_0\in\rA_*^+$ such that $\sup\limits_{\ph\in
K}\vert\ph(a)\vert\to 0$ when $\psi_0(aa^*+a^*a)\to 0$, $a\in
B_\srA$.

\proof The new ingredient will be the proof of (iii) $\implies$
(iv); the proofs of the other implications are standard.

(i) $\implies$ (ii): Let $(p_n)$ be a disjoint sequence of
projections in $\rA$: they generate an abelian von Neumann
subalgebra $\B$ of $\rA$. Let $\rho$ be the restriction map
$\rA_*\to\B_*$, $\ph\mapsto {\ph_\mid}_\B$. Then $\rho(K)$ is
weakly relatively compact in $\B_*\sim\ell_1$, and consequently
(by the commutative result, i. e. Dunford-Pettis Theorem, see [D], p. 93),
$\sup\{\vert f(p_n)\vert\,\mid f\in
\rho(K)\,\}\to 0$, i.e. $\sup\{\vert \ph(p_n)\vert\,\mid \ph\in
K\,\}\to 0$.

(ii) $\implies$ (iii): Assume that for some disjoint sequence
$(p_n)$ of projections in $\rA$ and some sequence $(\ph_n)$ in $K$
we have $\inf\limits_n\N{p_n\ph_n p_n}=\delta>0$. For every $n$ we
have $\N{p_n\ph_n p_n}=\langle p_n\ph_n p_n,u_n\rangle=
\ph_n(p_nu_np_n)$ for some partial isometry $u_n$ in $\rA$. Let
$p_nu_np_n=a_n+ib_n$ be the decomposition into  real and imaginary
parts. Then $a_n,b_n$ belong to  the unit ball of $p_n\rA p_n$.
Thus one of the sets $N_1=\{n\geq 1\mid
\vert\ph_n(a_n)\vert\geq\delta/2\}$ and $N_2=\{n\geq 1\mid
\vert\ph_n(b_n)\vert\geq\delta/2\}$ is infinite. Suppose w.l.o.g.
that $\vert\ph_n(a_n)\vert\geq\delta/2$ for every $n\geq 1$. The
von Neumann algebra $\C$ generated by the $a_n$'s is commutative
(since they are hermitian and disjoint). The image $\rho(K)$ of
the set  $K$ by the restriction map $\rho: \rA_*\to\C_*$ still
verifies (ii). However, in the commutative case, it is clear that
 (ii) implies (iii). Therefore, $\langle
\ph_n,a_n\rangle= \langle \rho(\ph_n),a_n\rangle\to 0$ when
$n\to\infty$, a contradiction.

(iii) $\implies$ (iv): if $K$ is 1-equiintegrable, then by Prop. 4.13
for every $n\geq 1$ there is $\ph_n\in\rA_*^+$ such that for every $\ph\in
K$ there exist  $x,y\in B_\srA$ such that $\N{\ph-(x\cdot
\ph_n+\ph_n\cdot y)}<1/n$. If $a\in B_\srA$ we have then:
 $$\eqalign{\abs{\ph(a)}
 &\leq\frac 1n+\vert \langle x\cdot\ph_n
 +\ph_n\cdot y,a\rangle\vert\leq \frac 1n
 +\vert\ph_n(ax)\vert+\vert\ph_n(ya)\vert\cr
 &\leq \frac 1n+\ph_n(aa^*)^{1/2}\ph_n(x^*x)^{1/2}
 +\ph_n(a^*a)^{1/2}\ph_n(yy^*)^{1/2}\cr
 &\leq \frac 1n+(2\|\ph_n\|)^{1/2}(\ph_n(aa^*)
 +\ph_n(a^*a))^{1/2}}$$
Set $\psi_0=\sum\limits_{n\geq 1}2^{-n}\N{\ph_n}^{-1}\ph_n$. Then
if
 $\psi_0(aa^*+a^*a)< 2^{-n-1}\N{\ph_n}^{-2}n^{-2}$,
we obtain $|\ph(a)|\leq 2/n$ for every $\ph\in K$, and thus prove
(iv).

(iv) $\implies$ (i): Let $f\in \rA^*$ be any w*-limit point of
$K$. Then $\abs{f(a)}\to 0$ when $\psi_0(aa^*+a^*a)\to 0$, $a\in
B_\srA$. Consequently, the linear functional $f$ is
strong*-continuous on bounded sets of $\rA$, so it is
w*-continuous ([T], Theorem II.2.6) and thus belongs to
$\rA_*$.\cqfd\medskip

\noindent{\bf Remark.} The equivalence of conditions (i) and (iii)
in Theorem 4.14 was also obtained in [HRS] for finite von Neumann
algebras.

\medskip
 Another application of Proposition 4.13  is the following well
known result of H. Jarchow [Ja]:

\proclaim Theorem 4.15. Every reflexive subspace of the predual of
a von Neumann algebra is super-reflexive.

\proof Since a Banach space is reflexive iff its unit ball is
weakly compact (or, equivalently, weakly relatively compact), then
by Theorem 4.14,  a closed subspace of a predual of von Neumann
algebra is reflexive iff its unit ball is $1$-equiintegrable.  Let
now $X$ be a reflexive subspace of the predual $\rA_*$. It is
super-reflexive iff all its ultrapowers are reflexive. Such an
ultrapower $X_\U$ is a closed subspace of $(\rA_*)_\U$ which we
identify with $\A_*$. Let $\eps>0$ and $\ph_\eps\in \rA_*$ such
that $d(\ph, B_\srA\cdot \ph_\eps+\ph_\eps\cdot B_\srA)<\eps$ for
every $\ph\in B_X$. Then clearly $d(\t \ph, B_{\srA_\U}\cdot
\hat\ph_\eps+\hat\ph_\eps\cdot B_{\srA_\U})\leq\eps$ for every $\t
\ph\in B_{X_\U}$, where $\hat \ph_\eps$ is the canonical image of
$\ph_\eps$ in $(\rA_*)_\U$. Thus $d(\t \ph, B_{\A}\cdot
\hat\ph_\eps+\ph_\eps\cdot B_{\A})\leq\eps$ for every $\t \ph\in
B_{X_\U}$, and so $B_{X_\U}$ is 1-equiintegrable, hence weakly
relatively compact.\cqfd

\noindent{\bf Remark.} It is easy to see that if $\psi_0$ is the
``control measure'' for $B_X$ given by Akemann's condition (iv) in
Theorem 4.14, then $\widehat\psi_0$ is a control measure for
$B_{X_\U}$ (in virtue of the strong*-density of $B_{\srA_\U}$ in
$B_\A$).

\bigskip\centerline {\bf 5. Subspaces containing $\ell_p$}

\medskip
The following is the main result of this section. It gives several
characterizations of the subspaces of $L_p(\rA)$ which contain
$\ell_p$.

\proclaim Theorem {5.1}. Let  $0<p<\infty$, $p\neq 2$ and $\rA$ be
a von Neumann algebra. Let $X\subset L_p(\rA)$ be a closed
subspace. The following statements are equivalent:\hfill\break
 i) $X$ contains an almost disjoint normalized sequence.\hfill\break
 ii) $X$ contains a basic sequence asymptotically $1$-equivalent to
the $\ell_p$-basis (and, if $1\leq p<\infty$,  spanning an almost
$1$-complemented subspace of $L_p(\rA)$).\hfill\break
 iii) $X$ contains a subspace isomorphic to $\ell_p$.\hfill\break
 iv) $X$ contains uniformly the spaces $\ell_p^n$, $n\geq 1$.\hfill\break
 v) For some $q\in(0,p)$, (or equivalently, for every $0<q<p$) and
for every $h\in L_r(\rA)$, where $\frac 1r=\frac 1q-\frac 1p$, the
restriction $T_{h,p,q}{\big|}_X$ is not an isomorphism,
where:\hfill\smallskip
 \centerline{$T_{h,p,q}: L_p(\rA)\To
L_q(\rA)\oplus L_q(\rA),\ x\mapsto (xh,hx)$} \hfill\break
 (If $\rA$ is $\sigma$-finite and $h_0$ is an element of $L_r(\rA)_+$
with full support, it is sufficient to test this condition on
$T_{h_0,p,q}$).\hfill\break
 vi) If in addition $0<p<2$: the unit ball of $X$ is
not $p$-equiintegrable.

To prove the last equivalence in Theorem 5.1, we shall need the
following lemma. As usual, we denote by $(\eps_n)$ a sequence of
independent Bernouilli variables (random signs) and by ${\bf
E}_\eps$ the corresponding expectation.

\proclaim Lemma {5.2}. Let $0<p<2$  and $K$ be a
$p$-equiintegrable sequence in $L_p(\rA)$. Then\hfill\break
\centerline{$ \lim\limits_{n\to\infty}n^{-1/p}\sup\,\{\,{\bf
E_\eps\,}\Vert \sum\limits_{i=1}^n \eps_i h_i\Vert\,\mid
h_1,...,h_n\in K\}=0$}

\proof For every $\alpha>0$, we can find by Proposition 4.13  an
element $h_0\in L_p(\rA)_+$, and for every $h\in K$, elements
$x,y$ in the unit ball of $\rA$ such that
$\N{h-xh_0-h_0y}<\alpha$. Given $h_1,...,h_n\in K$ let
$x_1,y_1,...,x_n,y_n\in B_\srA$ such that
$\N{h_i-x_ih_0-h_0y_i}<\alpha$, $i=1,...,n$. Let $r>0$ such that
$1/p=1/2+1/r$. We have:
 $${\bf E}_\eps\,\Vert\sum_{i=1}^n\eps_ix_ih_0\Vert_p
 \leq {\bf E}_\eps\,\Vert\sum_{i=1}^n\eps_ix_ih_0^{p/2}\Vert_2
 \Vert h_0^{p/r}\Vert_r
 =\Vert h_0\Vert_p^{p/r}(\sum_{i=1}^n
 \Vert x_ih_0^{p/2}\Vert^2_2)^{1/2}
 \leq n^{1/2}\N{h_0}_p$$
and similarly,
 $${\bf E}_\eps\,\Vert\sum\limits_{i=1}^n\eps_ih_0y_i\Vert_p
 \leq n^{1/2}\N{h_0}_p$$
Since $L_p(\rA)$ is of type $p$ with constant $1$ (this follows
from interpolation if $1<p<2$ and from the $p$-norm inequality if
$0<p<1$), we have
 $${\bf E}_\eps\,\Vert\sum_{i=1}^n\eps_i(h_i-x_ih_0
 -h_0y_i)\Vert_p
 \leq(\sum_{i=1}^n\Vert h_i-x_ih_0
 -h_0y_i\Vert_p^p)^{1/p}\leq n^{1/p}\alpha$$
Combining the preceding inequalities, we deduce
 $$\limsup\limits_{n\to\infty} n^{-1/p}
 \sup\,\{\,{\bf E_\eps\,}\Vert\sum\limits_{i=1}^n
 \eps_i h_i\Vert\,\mid h_1,...,h_n\in K\}\leq\alpha$$
which proves the lemma. \cqfd

\bg\noindent{\bf Proof of Theorem 5.1.} The implication (i)
$\implies$ (ii) is the implication (iii)$\implies$ (iv) in the Prop. 2.11;
the implications (ii) $\implies$ (iii) $\implies$ (iv) are
trivial. The implication (iv) $\implies$ (i) is a special case of
Theorem 3.1.

(i) $\implies$ (v): if $(a_n)$ is a normalized disjoint sequence,
$a_nh\to 0$ and $ha_n\to 0$ for every $h\in L_r$. For $s_nh\to 0$
(resp. $hs_n\to 0$) in $L_r$ for every disjoint sequence $(s_n)$
of projections (since the set $\{h\}$ is $p$-biequiintegrable).

(v) $\implies$ (i): let $\cal F$ be the set of finite subsets of
$L_r(\rA)$, ordered by inclusion. Let $\Phi$ be the set of cofinal
subsets of $\cal F$ and  $\U$  an ultrafilter on $\cal F$
containing $\Phi$ (see the proof of Proposition 4.13). Note that
the ultrafilter $\U$ is necessarily countably incomplete: so we
can find a family $(\eps_F)$ of strictly positive real numbers
such that $\lim\limits_{F,\U}\eps_F=0$. By (v) we can choose for
every $F\in\cal F$ an element $x_F\in X$, with $\N{x_F}_p=1$, such
that $\N{x_Fh}<\eps_F$ and $\N{hx_F}<\eps_F$ for every $h\in F$
(apply hypothesis (v)  to $h_F=(\sum\limits_{h\in
F}(\abs{h}^2+\vert h^*\vert^2))^{1/2}$ ). It follows that the
family $(x_F)$ defines an element $\xi$ of $X_\U$ (hence of
$L_p(\rA)_\U$) which verifies $\hat h\xi=0=\xi\hat h$ for every
$h\in L_r(\rA)$. Consequently, $\xi$ is disjoint from $L_p(\rA)$,
and so by Theorem 2.7 we can extract from the family $(x_F)$ an
almost disjoint sequence. Thus we get (i). In the case where $\rA$
is $\sigma$-finite and $h_0\in L_r(\rA)_+$ with full support such
that $T_{h_0,p,q}$ is not an isomorphism, we choose a normalized
sequence $x_n$ in $X$ such that $\N{x_nh_0}\to 0$ and
$\N{h_0x_n}\to 0$, and consider the element $\xi$ defined by the
sequence $(x_n)$ in some ultrapower $X_\U$ (associated with a free
ultrafilter over $\nat$). Then $\xi$ is disjoint from $h_0$ and
consequently from $L_r(\rA)$ (since $h_0$ has full support), and
finally from $L_p(\rA)$.

(i) $\implies$ (vi): This is clear since a disjoint normalized
sequence is not $p$-equiintegrable (nor is an almost disjoint
normalized sequence).

(vi) $\implies$ (i): if the unit ball of $X$ is not
$p$-equiintegrable, we can find a sequence $(h_n)$ of normalized
elements of $X$ and a  disjoint sequence $(p_n)$ of projections of
$\rA$ such that $\N{p_nh_np_n}_p>\delta>0$ for every $n\geq 1$. By
the Subsequence Splitting Lemma, we may suppose that
$h_n=h'_n+h''_n$, where $(h'_n)$ is $p$-equiintegrable and
$(h''_n)$ is disjoint. We have $p_nh'_np_n\to 0$, so we may
suppose that $\N{p_nh''_np_n}>\delta$, and consequently
$\N{h''_n}>\delta$ for every $n\geq 1$. Using  Lemma 5.2, we can
construct inductively a sequence $I_1<...<I_n<...$ of disjoint
intervals of $\nat$ and a sequence of signs $(\eps_i)$ such that
$\abs{I_n}^{-1/p}\Vert\sum\limits_{i\in I_n}
\eps_ih'_i\Vert<2^{-n}$ for every $n\geq 1$. Let
 $$a'_n=\abs{I_n}^{-1/p}\sum\limits_{i\in I_n}\eps_ih'_i,\quad
 a''_n=\abs{I_n}^{-1/p}\sum\limits_{i\in I_n} \eps_ih''_i\quad
 \hbox{ and}\quad a_n=a'_n+a''_n$$
Then $(a_n)\subset X$, $(a''_n)$ is equivalent to the
$\ell_p$-basis, and by a standard perturbation argument,
$(a_n)_{n\geq n_0}$ is also equivalent to the $\ell_p$-basis for
sufficiently large $n_0$. \cqfd

\smallskip
The equivalence between (i) and (vi) in Theorem 5.1 can be
extended to sequences in the following way.

\proclaim Proposition {5.3}. Let $0<p<2$, and $(h_n)\subset
L_p(\rA)$ be a bounded sequence. If $1<p<2$, suppose in addition
that $(h_n)$ is unconditional. Then the following assertions are
equivalent:\hfill\break i) $(h_n)$ is not
$p$-equiintegrable;\hfill\break ii) $(h_n)$ contains a subsequence
equivalent to the $\ell_p$-basis.

\proof That (ii) $\implies$ (i) is a consequence of Lemma 5.2. The
converse implication can be proved using arguments similar to
those used in the semi-finite case by [HRS] for the case $1\leq
p<2$ or [SX] for the case $p\leq 1$. We sketch these arguments for
the convenience of the reader (with a modified, somewhat shortened
proof in the case $0<p<1$). Assuming (i), we can choose, by
Corollary 4.9,  a subsequence of $(h_n)$ (for simplicity of
notation, we shall assume that it is $(h_n)$ itself), and a
 disjoint sequence of projections $(e_n)$ such that
$\N{e_nh_ne_n}>\delta>0$ for every $n\geq 1$.

\noindent a) The case $1< p< 2$.

For every finite sequence $(\lambda_n)$ of scalars, since the
projections $e_j$ are pairwise disjoint and $p\geq 1$, we have:
 $$\eqalign{{\bf E}_\eps\Vert \sum_n\eps_n\lambda_n h_n\Vert^p
 &\geq {\bf E}_\eps\sum_{j=1}^\infty\Vert e_j
 (\sum_n\eps_n\lambda_nh_n)e_j\Vert^p
 =\sum_{j=1}^\infty{\bf E}_\eps\Vert e_j(\sum_n\eps_n
 \lambda_nh_n)e_j\Vert^p\cr
 &\ge \sum_{j=1}^\infty\Vert \lambda_je_jh_je_j\Vert^p
 \geq \delta^p\sum_j\abs{\lambda_j}^p}$$
Thus by the unconditionality of the sequence $(h_n)$, we deduce
that $\Vert \sum\limits_n\lambda_n h_n\Vert^p\geq
c\sum\limits_n\abs{\lambda_n}^p$ for some $c>0$. The converse
inequality follows from the type $p$ property of $L_p(\rA)$.

\noindent b) The case $p\leq 1$.

Let $\U$ be a free ultrafilter over $\nat$; let $\t h$ be the
element of $L_p(\rA)_\U=L_p(\A)$ represented by the sequence
$(h_n)$, and for each $m$,  let $\hat e_m$ be the canonical image
of $e_m$ in $\rA_\U$. We have $\lim\limits_{m\to \infty} \Vert\hat
e_m\t h\Vert=0$, since the projections $\hat e_m$ are pairwise
disjoint. Similarly $\lim\limits_{m\to \infty} \Vert\t h\hat
e_m\Vert=0$. Let $m_k$ be such that $\Vert\hat e_m\t
h\Vert+\Vert\t h\hat e_m\Vert<\eps\cdot 2^{-k-1}$ for every $m\geq
m_k$. On the other hand, $\lim\limits_{n\to \infty} \Vert e_n
h_m\Vert=0=\lim\limits_{n\to\infty} \Vert h_m e_n \Vert$ for every
$m\in \nat$. So we can define inductively an increasing sequence
$(n_k)$, with $n_k\geq m_k$ for all $k$, in the following way:
$n_1=m_1$, and $n_{k+1}\geq\max(n_k+1,m_{k+1})$ is chosen such
that for every $j=1,...,k$:
 $$\eqalign{
 &i)\ \max(\Vert e_{n_{k+1}}h_{n_j}\Vert,
 \Vert h_{n_j}e_{n_{k+1}}\Vert)\leq \eps 2^{-k-1}\cr
 &ii)\ \max(\Vert h_{n_{k+1}}e_{n_j}\Vert,
 \Vert e_{n_j} h_{n_{k+1}}\Vert)\leq \eps2^{-j}}$$
Then $\max(\Vert e_{n_{j}}h_{n_k}\Vert,\Vert
h_{n_k}e_{n_{j}}\Vert)\leq \eps2^{-j}$ for every $j\neq k$. Set
$e=\sum\limits_{j\geq 1} e_{n_j}$ and $K=\sup\Vert h_n\Vert$. We
have for every $(\lambda_j)\in\ell_p$:
 $$\eqalign{K^p\sum_{j\geq 1}\vert \lambda_j\vert^p
 &\geq \Vert \sum_{j\geq 1}\lambda_j h_{n_j}\Vert^p
 \geq \Vert e(\sum_{j\geq 1}\lambda_jh_{n_j})e\Vert^p\cr
 &\geq \Vert \sum_{j\geq 1}\lambda_je_{n_j}h_{n_j}e_{n_j}\Vert^p
 -\sum_{j\geq 1}\sum_{k\neq j}\Vert\lambda_j e_{n_k}h_{n_j}e\Vert^p
 -\sum_{j\geq 1}\sum_{k\neq j}\Vert\lambda_j e_{n_j}h_{n_j}e_{n_k}\Vert^p\cr
 &\geq(\delta^p-2C^p\eps^p)\sum_{j\geq 1}\vert \lambda_j\vert^p}$$
where $C^p=\sum\limits_k 2^{-kp}$. Hence if
$\eps<2^{-1/p}(\delta/C)$, the sequence $(h_{n_j})$ is equivalent
to the $\ell_p$-basis. \cqfd

\medskip

\noindent{\bf Application: Kadec-Pe\l czy\'nski dichotomy for
non-commutative $L_p$-spaces}

\medskip

\noindent {\bf Proof of Theorem 0.2.} Theorem 0.2 follows from the
special case $q=2$ of Theorem 5.1. Note that if for some $h\in
L_r(\rA)$ the map $T_{h,p,2}: L_p(\rA)\To L_2(\rA)\oplus_2
L_2(\rA)$ restricts to an isomorphism $S=T_{h,p,2}\mid_X$, then
$X$ is isomorphic to the Hilbert space $S(X)$; if
$P:L_2(\rA)\oplus_2 L_2(\rA)\to S(X)$ is the orthogonal
projection, then $Q:=S^{-1}PT_{h,p,2}$ is a projection from
$L_p(\rA)$ onto $X$.\cqfd

\smallskip
 Another version of the Kadec-Pe\l czy\'nski dichotomy,
which is well known in the commutative case,  deals with
unconditional sequences.

\proclaim Proposition {5.4}. Let $2<p<\infty$ and $\rA$ be a von
Neumann algebra. Then every semi-normalized unconditional sequence
in $L_p(\rA)$ either is equivalent to the $\ell_2$-basis or has a
subsequence which is asymptotically $1$-equivalent to the
$\ell_p$-basis and spans a complemented subspace.

\proof Let $(h_n)$ be a semi-normalized unconditional sequence in
$L_p(\rA)$ (by semi-normalized we mean that
$0<\delta=\inf\limits_n\Vert h_n\Vert_p\le
M=\sup\N{h_n}_p<\infty$). Let $T_{h,p,2}: L_p(\rA)\To
L_2(\rA)\oplus_2 L_2(\rA)$ be as before. If for some $h\in
L_r(\rA)$ (with $1/r=1/2-1/p$), $\inf\{\Vert
T_{h,p,2}h_n\Vert_2\mid n\geq 1\}=c>0$, then $(h_n)$ is equivalent
to the $\ell_2$-basis. Indeed, then for every finite sequence
$(\lambda_n)$ of scalars:
 $$\eqalign{\sqrt 2\N{h}_r{\bf E}_\eps\Vert\sum_n\lambda_n\eps_nh_n\Vert_p
 &\geq {\bf
 E}_\eps\Vert\sum_n\lambda_n\eps_nT_{h,p,2}h_n\Vert_2\cr
 &=(\sum_n\vert\lambda_n\vert^2\Vert T_{h,p,2}h_n\Vert^2)^{1/2}
 \geq c\,(\sum_n\vert\lambda_n\vert^2)^{1/2}}$$
on the other hand,  by the type 2 property of $L_p(\rA)$,
 $${\bf E}_\eps\Vert\sum_n\lambda_n\eps_nh_n\Vert_p
 \leq C(\sum_n\vert\lambda_n\vert^2\Vert h_n\Vert^2)^{1/2}
 \leq CM(\sum_n\vert\lambda_n\vert^2)^{1/2}$$
If at the contrary we have $\inf\{\Vert T_{h,p,2}h_n\Vert_2\mid
n\geq 1\}=0$ for any $h\in L_r(\rA)$, then we can adapt
the proof of (v)$\implies$ (i) of Theorem 5.1  by choosing the
$x_F$'s in the sequence $(h_n)$. Then we deduce that $(h_n)$ has
an almost disjoint subsequence.\cqfd

\medskip

\noindent {\bf Remark.} The results of [HRS] concerning the
Banach-Saks properties of non-commutative $L_p$-spaces associated
with a finite von Neumann algebra can be extended to the present
 setting with the same proof (using our Lemma 5.2 in place
of Lemma 3.1 of [HRS]). We refer to [HRS], Definition 5.5 for the
definition of the various Banach-Saks properties (the terminology
is not completely fixed and their Banach-Saks property is
sometimes called ``weak Banach-Saks'', or
``Banach-Saks-Rosenthal'' property). Then $L_1(\rA)$ has the
Banach-Saks property, $L_p(\rA)$ has the $p$-Banach-Saks property
if $1\leq p\leq 2$ and the $2$-Banach-Saks property if $2\leq
p<\infty$; any $p$-equiintegrable weakly null sequence of
$L_p(\rA)$, $1<p<2$, has a strong  $p$-Banach-Saks subsequence,
and a closed linear subspace of $L_p(\rA)$, $1<p<2$ has the strong
$p$-Banach-Saks property iff it has no subspace isomorphic to
$\ell_p$.

\medskip

\noindent{\bf Historical comments.}
 i) The commutative forerunner of Theorem 5.1 is due to H. P. Rosenthal
[Ro]. For finite von Neumann algebras,  Theorem 5.1 was proved in [HRS]
for $1\leq p<2$ and in [SX] for $p<1$. The proofs in both papers use
the notion of the $p$-modulus of uniform integrability, the
definition of which involves the trace. An analogue of this
modulus could be defined in the non-tracial, $\sigma$-finite case
too, via a normal faithful state, but it would have less tractable
properties (in particular with respect to conjugation and absolute
value). For finite von Neumann algebras,  some equivalences in
Theorem 5.1 were obtained  in [R1-2].\sn
 ii) Lemma 5.2 and Proposition 5.3 above are simply extensions
to our context of the corresponding results for finite von Neumann
algebras in [HRS]. \sn
 iii) Theorem 0.2 and Proposition 5.4 were proved in the commutative
context in the well-known paper of M. I. Kadec and A. Pe\l czy\'nski
[KaP]. These results were proved by  Sukochev [S]  in the case of
a finite von Neumann algebra, and by N. Randrianantoanina [R2] in
the semi-finite case (even in the more general setting of spaces
$E(\rA,\tau)$ associated with an order-continuous type 2 r.i.
function space $E$).

\bigskip\centerline{\bf 6. Operator space version}

\medskip
  This section is devoted to the analogues of
Theorems 3.1 and 0.2 in the category of operator spaces. Our
references for Operator Space Theory are [ER] and [P2]. Recall
that an operator space $E$ is a closed subspace of $B(H)$ for some
Hilbert space $H$, equipped with a natural sequence of matrix
norms. Let $M_n$ denote the space of complex $n\times n$ matrices
and $M_n(E)=M_n\otimes E$ the space of $n\times n$ matrices with
entries in $E$. As usual, $M_n(B(H))$ is identified with
$B(\ell_2^n(H))$ and the matrix norm on $M_n(E)$ is the one
induced by the natural inclusion of $M_n(E)$ into $M_n(B(H))$. An
abstract characterization of operator spaces was given by Ruan
(see [ER2]): a Banach space $E$ equipped with a sequence $(\N{\
}_n)$ of norms on the $M_n(E)$ can be identified with an operator
space iff the matricial norms $(\N{\ }_n)$ satisfy two  simple
conditions (``Ruan's axioms'').

Now let $E,F$ be two operator spaces;  a linear map $T:E\to F$ is
said to be completely bounded if
$$\N{T}_{cb}:=\sup_{n\geq 1}\N{{\rm id}_{M_n}\otimes T: M_n(E)\to
M_n(F)}<\infty $$ Then $\N{\ }_{cb} $ defines a norm on the space
$CB(E,F)$ of completely bounded operators from $E$ into $F$. An
operator $T:E\to F$ is called a complete isomorphism if it is a
linear bijection such that $T$ and $T^{-1}$ are completely
bounded. Two operator spaces $E$, $F$ are said to be completely
isomorphic (resp. $K$-completely isomorphic) if there exists a
complete isomorphism $T: E\to F$  (resp. with $\N{T}_{cb}\Vert
T^{-1}\Vert_{cb}\leq K$). Similarly, a linear subspace $F$ of $E$
is said to be completely complemented (resp. $K$-completely
complemented) in $F$ if there is a completely bounded projection
$P: E\to F$ (resp. with $\N{P}_{cb}\leq K$). 

Now let $\rA$ be a von Neumann algebra and $1\leq p\leq\infty$. We
will consider the natural operator space structure on $L_p(\rA)$
as introduced in [P1-2]. For $p=\infty$ a realization of $\rA$ as a
concrete von Neumann algebra, i.e. a unital w*-closed *-subalgebra of
some $B(H)$, gives an operator space structure on $\rA=L_\infty(\rA)$
(independent of  the realization since *-isomorphisms are completely
isometric). A standard operator space structure on the dual space $\rA^*$
follows by Operator Spaces Theory ($M_n(\rA^*)$ is identified with
the space $CB(\rA,M_n)$ and the corresponding sequence of
matricial norms satisfies Ruan's axioms). A specific  operator
space structure on $L_1(\rA)=\rA_*$ is induced by the natural
embedding of $\rA_*$ into its bidual $\rA^*$. In fact, as explained
in [P2], \S7 it is more convenient to consider $L_1(\rA)$ as the predual
of the opposite von Neumann algebra $\rA^{op}$, which is isometric
(but not completely isomorphic) to $\rA$, and to equip  $L_1(\rA)$
with the operator space structure inherited from $(\rA^{op})^*$.
The main reason for this choice is that it insures that the equality
$L_1(M_n\otimes \rA)=S_1^n\widehat \otimes L_1(\rA)$ (operator
space projective tensor product) holds true (see [Ju3], \S3). Finally the
operator space structure of $L_p(\rA)$ is obtained by complex
interpolation, using the well known interpretation of $L_p(\rA)$ as
interpolation space $(\rA,L_1(\rA))_{1/p}$ (see [Te2]). 

We will need the following convenient characterization of the
operator space structure of the subspaces of $L_p(\rA)$. Note that there
is a natural {\it algebraic} identification of $L_p(M_n\otimes\rA)$
with $M_n(L_p(\rA))$. Following [P1], if $E$ is an operator space one
sets $S_p^n[E]:=(S_\infty^n[E],S_1^n[E])_{1/p}=(M_n[E],S_1^n\widehat
\otimes E)_{1/p}$; by [P1], Cor 1.4 we have when $(E_0,E_1)$ is a
compatible interpolation couple:
$S_p^n[(E_0,E_1)_{1/p}]=(M_n[E_0],S_1^n\widehat
\otimes E_1)_{1/p}$. Consequently we have:
$S_p^n[L_p(\rA)]=(M_n\otimes\rA,L_1(M_n\otimes\rA))_{1/p}=L_p(M_n\otimes
\rA)$; in other words $S_p^n[\rA]$ identifies with  the linear space
$M_n(L_p(\rA))$ equipped with the norm of $L_p(M_n\otimes \rA)$. Note
that if $\rA$ is commutative, say
$\rA=L_p(\Omega,\mu)$ for some measure space $(\Omega,\mu)$, it
turns out that $S_p^n[L_p(\rA)]=L_p(\Omega,\mu;S_p^n)$, the space
of $p$-integrable functions with values in $S_p^n$. Recall also that if
$F$ is a closed linear subspace of $E$, the norm on $S_p^n[F]$ is
induced from that of $S_p^n[E]$.  The norms on
$S_p^n[E],n\geq 1$, completely determine the operator space structure of
$E$ in the following sense (see [P1], prop. 2.3):

\proclaim Lemma {6.1}. Let $E_1$ and $E_2$ be two operator spaces. Then a
linear map $T: E_1\to E_2$ is completely bounded iff \hfill\break
\centerline{$\sup\limits_n\Vert{\rm id}_{S_p^n}\otimes T:
S_p^n[E_1]\to S_p^n[E_2]\Vert<\infty$} moreover in this case the
supremum above is equal to $\N{T}_{cb}$.

The embedding results in sections 3 and 5 can be improved into
results for the category of operator spaces. We first consider
subspaces of $L_p(\rA)$ containing $\ell_p$.

\proclaim Theorem {6.2}. Let $0<p<\infty$, $p\neq 2$ and $X$ be a
closed subspace of $L_p(\rA)$. If $X$ contains uniformly the
spaces $\ell_p^n$, $n\geq 1$ as Banach spaces, then given any
$\eps>0$, $X$ contains a subspace $(1+\eps)$-completely isomorphic
to $\ell_p$ and $(1+\eps)$-completely complemented in $L_p(\rA)$.

As a corollary we immediately get the following operator space
version of the Kadec-Pe\l czy\'nski dichotomy:

\proclaim Corollary {6.3}. Let $2<p<\infty$, and $X$ be a closed
subspace of $L_p(\rA)$. Then either $X$ is (Banach) isomorphic to
a Hilbert space or $X$ contains a subspace which is completely
isomorphic to $\ell_p$ and completely complemented in $L_p(\rA)$.

\noindent{\bf Remark.} Corollary 6.3. has been known to M. Junge
and the second author for a semifinite $\rA$ (and also when $\rA$
is a type III algebra of some particular form).\bigskip

Now we turn to the operator space analogue of Theorem 3.1. In the
following, given an operator space $F$, the spaces $\ell_p(F)$ and
$\ell_p^n(F)$ are equipped with the natural operator space
structure introduced in [P1] (via complex interpolation). More
generally, if $(F_j)_{j\geq 1}$ is a sequence of operator spaces,
the space $(\bigoplus\limits_{j\geq 1} F_j)_p$ has also a natural
operator space structure.

\proclaim Theorem {6.4}. Let $1\leq p<\infty$, and $X$ be a closed
subspace of $L_p(\rA)$. Let $(F_j)_{j\geq 1}$ be a sequence of
finite dimensional operator spaces. Assume that there is a
constant $K$ such that for all $n,j\geq 1$, $X$ contains a
subspace $Y_{j,n}$which is $K$-completely isomorphic to
$\ell_p^n(F_j)$. \hfill\break
 i) Then for every $\eps>0$, $X$
contains a subspace $(K+\eps)$-completely isomorphic to
$F=(\bigoplus\limits_{j\geq 1} F_j)_p\,$. \hfill\break
 ii) If in addition each $Y_{j,n}$ is $C$-completely complemented in
$L_p(\rA)$ then $X$ contains a subspace $(K+\eps)$-completely
isomorphic to $F$ and $(CK+\eps)$-completely  complemented in
$L_p(\rA)$.

Specializing to  Schatten classes we get the following:

\proclaim Corollary {6.5}. Let $1\leq p<\infty$, and $X$ be a
closed subspace of $L_p(\rA)$. If $X$ contains subspaces uniformly
completely isomorphic to $S_p^n$, $n\geq 1$, (resp. and uniformly
completely complemented in $L_p(\rA)$) then $X$ contains a
subspace completely isomorphic to $K_p=(\bigoplus\limits_{n\geq 1}
S_p^n)_p$ (resp. and completely complemented in $L_p(\rA)$).

\smallskip\noindent{\bf Remarks.} i) Corollary 6.5 can be used to
simplify some proofs in [JNRS].

ii)  It is worth noting that contrary to Theorem 6.4, the
assumption in Theorem 6.2 is only at the Banach space level! (So
the latter cannot be considered as a special case of the former.)

\smallskip

The proofs of Theorems 6.2 and 6.4 are very similar and that of
Theorem 6.2 is simpler, so we give only the proof of the
latter.\medskip

\noindent{\bf Proof of Theorem 6.4.} By the proof of Theorem 3.1,
given a sequence $(\eps_j)$ of positive real numbers (the
$\eps_j$'s being very small), there is a sequence $(E_j)$ of
finite dimensional subspaces of $X$ and a disjoint sequence
$(p_j)$ of projections of $\rA$, such that $E_j$ is $K$-isomorphic
to $F_j$ by some isomorphism $T_j: F_j\to E_j$, and such that
 $$\forall j\geq 1,\forall h\in E_j,\qquad
  \N{h-p_jhp_j}\leq\eps_j\N{h}\eqno (\ast)$$
Define $T: F\to E=\sum\limits_j E_j$, $x=(x_j)\mapsto
Tx=\sum\limits_j T_jx_j$. Then $T$ is an isomorphism, see the
proof of Theorem 3.1.

Reexamining that proof, we see that, under the present hypothesis,
each of the spaces $E_j$ constructed there is in fact
$K$-completely isomorphic to the corresponding $F_j$. More
precisely, $T_j$ can be defined so that
 $$\N{T_j^{-1}}_{cb}\leq 1\hbox{ and } \N{T_j}_{cb}
 \leq K$$
Indeed, keeping the notations used in the proofs of Lemma 3.2 and
Theorem 3.1, we have that $T_j$ is some $S_{j,i,n}$ at the
beginning of the proof of Theorem 3.1 and $E_j=S_{j,i,n}(F_j)$.
However, $S_{j,i,n}=T_{j,n}\,I_i$, where $I_i: F_j \to
\ell_p^n(F_j)$ is the natural embedding of $F_j$ into the i-th
coordinate of $\ell_p^n(F_j)$, i.e. $I_i(x)=e_i\otimes x$, and
where $T_{j,n}: \ell_p^n(F_j)\to X$ is the embedding  given by
Lemma 3.2. More precisely, by the discussion at the end of the
proof of Lemma 3.2, for every $n,\; j$ there are an integer $N$
and a linear map $S_n: \ell_p^n\to \ell_p^N$ which satisfy the
following: firstly, $S_n$ sends the basis vectors of $\ell_p^n$
into disjoint blocks of $\ell_p^N$; secondly,
$T_{j,n}=T^{j,N}\,\big(S_n\otimes {\rm Id}_{F_j}\big)$, where
$T^{j,N}:\ell_p^N(F_j)\to X$ is a $K$-complete embedding whose
existence is guaranteed by the assumption of Theorem 6.4.
Therefore, we obtain that
 $$S_{j,i,n}=T^{j,N}\,\big(S_n\otimes {\rm
 Id}_{F_j}\big)\,I_i\,.$$
Since both $ I_i$ and $S_n\otimes {\rm Id}_{F_j}$ are completely
isometric embeddings, we deduce the desired assertion on $T_j$.

We shall show that $T$ is now a complete isomorphism. To this end,
by Lemma 6.1 we need only to consider
 $${\rm id}_{S_p^m}\otimes T: S_p^m[F]\to S_p^m[E],
 \qquad m\geq 1$$
Fix $m\geq 1$ and let $\t p_j={\rm id}_{\ell_2^m}\otimes p_j$.
Then $(\t p_j)$ is a disjoint sequence of projections in
$M_n\otimes\rA$. We claim that (with $d_j=\dim E_j$):
 $$\forall m,j\geq 1,\forall h\in S_p^m[E_j],\qquad
 \N{h-\t p_jh\t p_j}_{S_p^m[E_j]}
 \leq d_j\eps_j\N{h}_{S_p^m[E_j]}\eqno(\dagger)$$
(Note that the constant on the right hand side does not depend on
$m$). Indeed, choose an Auerbach basis $(\xi_i)_{1\leq i\leq d_j}$
in $E_j$ and let $(\xi^*_i)_{1\leq i\leq d_m}$ be the dual basis
in the dual space $E_j^*$. For all $i,k=1,...,d_j$, $i\neq k$ we
have $\langle \xi_i,\xi_i^*\rangle=\N{\xi_i}=\Vert\xi_i^*\Vert=1$,
and $\langle \xi_i,\xi_k^*\rangle=0$. Every $h\in S_p^m[E_j]$ can
be written as $h=\sum\limits_{i=1}^{d_j} a_i\otimes \xi_i$, where
the $a_i$, $i=1,...,d_j$ belong to $ S_p^m$. Thus by $(\ast)$:
 $$\qquad\N{h-\t p_jh\t p_j}_{S_p^m[E_j]}
 \leq d_j\eps_j\sup_i\N{a_i}_{S_p^m}$$
Recall that any bounded functional on an operator space is
automatically completely bounded, and that its cb-norm is equal to
its norm. Hence
 $$\Vert{\rm id}_{S_p^m}\otimes \xi_i^*:
 S_p^m[E_j]\to S_p^m\Vert=1,\quad 1\leq i\leq d_j$$
whence
 $$\sup_i\N{a_i}_{S_p^m}\leq \N{h}_{S_p^m[E_j]}$$
Combining the previous inequalities we obtain our claim
$(\dagger)$. Now using $(\dagger)$ instead of $(\ast)$ and
repeating the arguments in the proof of Theorem 3.1 with ${\rm
id}_{S_p^m}\otimes T$ in place of $T$, we deduce that, if
$d_j\eps_j<1$ for each $j$ and $\ds\sum_j {d_j\eps_j\over
1-d_j\eps_j}=\eps<1$,
 $$\Vert {\rm id}_{S_p^m}\otimes T^{-1}\Vert
 \leq (1-\eps)^{-2},\qquad
 \Vert{\rm id}_{S_p^m}\otimes T\Vert
 \leq (1+\eps)K$$
Since $m\geq 1$ is arbitrary, we obtain:
 $$\Vert T^{-1}\Vert_{cb}\leq (1-\eps)^{-2}, 
 \qquad\N{T}_{cb}\leq(1+\eps)K$$
This proves the part (i) of Theorem 6.4. The part (ii) can be
proved similarly by combining the previous arguments with the
proof of the part iii) of Theorem 3.1.\cqfd

\bigskip\centerline{\bf Appendix:  Equality case in non-commutative
Clarkson inequality}

\medskip

 \proclaim Theorem A1. Let $\rA$ be a von Neumann algebra
and $0<p<\infty$, $p\neq 2$. Two elements $x,y$ of $L_p(\rA)$
verify the equality:\smallskip
\centerline{$\N{x+y}^p+\N{x-y}^p=2(\N{x}^p+\N{y}^p)$} \noindent if
and only if they are disjoint.

This result was stated by H. Kosaki [Ko2] in the case $p>2$, and
proved by reduction to the equality case of an inequality valid in
$L_{p/2}(\rA)_+$. We shall follow the same pattern, but the
argument is different when $p<2$. The equality case of this
auxiliary inequality is given by the following:

\proclaim Proposition A2. Let $\rA$ be a von Neumann algebra and
$0<r<\infty$, $r\neq 1$. Two positive elements $a,b$ of $L_r(\rA)$
verify the equality:\smallskip
\centerline{$\Tr(a+b)^r=\Tr(a^r)+\Tr(b^r)$} \noindent if and only
if they are disjoint.

We first deduce  Theorem A1 from Proposition A2. Let $r=p/2$ and
$a=x^*x$, $b=y^*y$. Then:
 $$\eqalign{\Tr (a^r)+\Tr(b^r)
 &= \N{x}^p+\N{y}^p=\frac 12(\N{x+y}^p+\N{x-y}^p)\cr
 & = \frac 12 \Tr[(a+b+(x^*y+y^*x))^r+(a+b-(x^*y+y^*x))^r]\cr
 &\left\{\eqalign{&\leq\Tr(a+b)^r \hbox{ if $0<r\leq 1$}\cr
                 & \geq \Tr(a+b)^r \hbox{ if $r>1$}}\right.
 }$$
where we have used the operator-concavity of the function
$t\mapsto t^r$ if $0<r\leq 1$ (see [B], chap. V), and the
convexity of the $L_r$-norm and of the function $t\mapsto t^r$ if
$r\geq 1$. Note that the reverse inequalities are always true:
 $$\Tr(a+b)^r\left\{\eqalign{&\leq \Tr (a^r)+\Tr(b^r)
  \hbox{ if $0<r\leq 1$}\cr
  &\geq \Tr (a^r)+\Tr(b^r)\hbox{ if $r>1$}}\right.$$
(see [Ko3] Lemma 3 in the first case and [Ko2] Lemma 3.3 in the
second case). So we are in the equality case
$$\Tr(a+b)^r=\Tr (a^r)+\Tr(b^r)$$
and by  Proposition A2 above this implies that $a$ and $b$ have disjoint
supports. Since the  support of $a$ (resp. b) coincides with the
right support of $x$ (resp of $y$) this means that $x$ and $y$
have disjoint right supports. Replacing $x,y$ by their conjugates,
we see that the left supports of $x$ and $y$ are disjoint too, so
$x\perp y$.\cqfd

\smallskip\noindent{\bf Proof of Proposition A2:} The case $r>1$
was treated in [Ko2] Proposition 6.3, using a differentiation
argument in the strictly convex Banach space $L_r(\rA)$. So we
consider only the case $0<r<1$.

>From $0\leq a,b\leq a+b$ we infer the existence of $c,d\in \rA$
with $\N{c}\leq 1$, $\N{d}\leq 1$ such that
 $$a^{1/2}=c(a+b)^{1/2},\quad b^{1/2}=d(a+b)^{1/2}$$
Hence
 $$\eqalign{a
 &=(a+b)^{1/2}c^*c(a+b)^{1/2}=c(a+b)c^*\cr
 b&=(a+b)^{1/2}d^*d(a+b)^{1/2}=d(a+b)d^*}$$
Note that we can choose $c, d$ such that $c^*c+d^*d=s(a+b)$. Since
$0\leq r\leq 1$ we have by Hansen's inequality (see [Han] for
bounded operators, and [Ko2] Lemma 3.6 for operators in
$L_p(\rA)$):
 $$\eqalign{
 &a^r=(c(a+b)c^*)^r\geq c(a+b)^rc^*\cr
 &b^r=(d(a+b)d^*)^r\geq d(a+b)^rd^*}$$
Hence
 $$\eqalign{\Tr(a^r+b^r)
 &\geq\Tr(c(a+b)^rc^*)+\Tr(d(a+b)^rd^*)\cr
 &=\Tr((a+b)^r(c^*c+d^*d))\cr
 &=\Tr(a+b)^r=\Tr(a^r+b^r)}$$
so the inequalities above become all equalities; then we have:
 $$\eqalign{&a^r= c(a+b)^rc^*\cr
 &b^r=d(a+b)^rd^*}$$
(since the differences are positive and of zero trace). We
distinguish now two cases:\medskip

\noindent \underbar{Case 1}: $r\leq 1/2$. Since $2r\leq 1$ we may
use Hansen's inequality again:
 $$a^r=c((a+b)^{1/2})^{2r}c^*\leq (c(a+b)^{1/2}c^*)^{2r}
 \leq (a^{1/2})^{2r}=a^r\eqno (\ast)$$
where the last inequality follows from the inequality
$c(a+b)^{1/2}c^*\leq (c(a+b)c^*)^{1/2}=a^{1/2}$ (Hansen's
inequality) and $2r\leq 1$. Therefore, the inequalities in
$(\ast)$ above are equalities:
 $$a^r=(c(a+b)^{1/2}c^*)^{2r}$$
whence
 $$a^{1/2}=c(a+b)^{1/2}c^*$$
equivalently:
 $$a^{1/2}=ca^{1/2}=a^{1/2}c^*$$
which implies in particular that $s(a)\leq r(c)$ and that
$a=cac^*$. Recalling that $a=c(a+b)c^*$, we see that $cbc^*=0$,
hence $r(c)\perp s(b)$. So finally $s(a)\perp s(b)$, which ends
the proof of case\nobreak~1.\medskip

\noindent \underbar{Case 2}: $1/2<r<1$. By the equalities $a^r=
c(a+b)^rc^*$ and $a^{1/2}=c(a+b)^{1/2}$ we have
 $$a^r=a^{1/2}(a+b)^{r-1/2}c^*$$
whence $a^{r-1/2}=s(a)(a+b)^{r-1/2}c^*$, and so:
 $$a^{2r-1}=c(a+b)^{r-1/2}s(a)(a+b)^{r-1/2}c^*
 \leq c(a+b)^{2r-1}c^*$$
but since $2r-1\leq 1$, we may use Hansen's inequality again:
 $$c(a+b)^{2r-1}c^*\leq (c(a+b)c^*)^{2r-1}=a^{2r-1}$$
and thus we obtain the equality:
 $$a^{2r-1}=c(a+b)^{2r-1}c^*$$
If $2r-1\leq 1/2$,  i.e. $r\leq 3/4$, then as in Case 1,  we
deduce that $a^{1/2}=c(a+b)^{1/2}c^*$ and then $s(a)\perp s(b)$.
If not, we iterate the procedure. Define the sequence $(r_n)$ by
$r_0=r$, $r_{n+1}=2r_n-1$. The interval $(1/2,1]$ contains
finitely many points of this sequence (which converges to
$-\infty$). Let $N$ be the first integer such that $r_N\leq 1/2$.
We have $0<r_N\leq 1/2$ and $1/2< r_n<1$ for $n=0,...,N-1$. So we
have inductively
 $$a^{r_n}=(c(a+b)c^*)^{r_n}$$
for $n=0,...,N$. For $n=N$ this equality implies
$a^{1/2}=c(a+b)^{1/2}c^*$ and finally that $s(a)\perp s(b)$. \cqfd

\bigskip

\bg\centerline{\eightsc References}\bigskip\nobreak

{\eightpoint \ref A&C. A. Akemann&The dual space of an operator
algebra,&Trans. Amer. Math. Soc. &{\bf 126} (1967) 268-302.

\ref AM&H. Araki, T. Masuda &Positive cones and $L^p$-spaces for
von Neumann algebras,&Publ. Res. Inst. Math.  Sci.&{\bf 18} (1982)
201-207.

\ref ArL&J. Arazy, J. Lindenstrauss&Some linear topological
properties of the space $C_p$ of operators on Hilbert
space,&Compositio Math.&{\bf 30} (1975) 81-111.

\ref B& R. Bhatia&& Matrix Analysis,&Graduate Texts in Mathematics
{\bf 169}, Springer-Verlag, New York,1997.

\ref BL& B. Beauzamy, J.-T. Laprest\'e&& Mod\`eles \'etal\'es des
espaces de Banach,&Travaux en Cours,
Hermann, Paris,1984.

\ref D& J. Diestel&& Sequences and Series in Banach Spaces,&Graduate Texts
in Mathematics {\bf 92}, Springer-Verlag, New York,1984.

\ref DJLT& P. N. Dowling, W. B. Johnson, C. J. Lennard, B. Turett& The
optimality of James's distortion theorem,&Proc. Amer. Math. Soc.& {\bf 125},
1997 (167-174).

\ref EJR& E. Effros, M. Junge, Z. J. Ruan&  Integral mappings and
the principle  of local reflexivity  for noncommutative
 $L^1$-spaces,& Ann. Math.& {\bf 151} (2000), 59-92.

\ref ER1& E. Effros,  Z.J. Ruan&${\cal O}{\cal L}_p$-spaces,&
Contemporary Math.&  {\bf  228} (1998) 51-77.

\ref ER2& E. Effros,  Z.J. Ruan&& Operator spaces,& London
Mathematical Society Monographs, New Series {\bf 23},  Oxford
University Press, New York, 2000.

\ref F& F. Fidaleo&Canonical operator space structures on
non-commutative $L_p$ spaces,& J. Funct. Anal.&{\bf 169} (1999),
226-250.

\ref G& U. Groh& Uniform ergodic theorems for identity preserving
Schwartz maps on W*-algebras,&J. Operator Theory& 11 (1984)
395-404.

\ref H& U. Haagerup& $L^{p}$-spaces associated with an arbitrary
von Neumann algebra,& Alg\`ebres d'op\'erateurs et leurs
applications en physique math\'ematique &({\sl Proc. Colloq.},
Marseille, 1977), pp. 175--184, {\sl Colloq. Internat. CNRS}, 274,
CNRS, Paris, 1979.

\ref HRS& U. Haagerup, H. Rosenthal, F.A. Sukochev& Banach
embedding properties of non-commutative $L^p$-spaces, & Memoirs
Amer. Math. Soc.& to appear.

\ref Han& F. Hansen&An operator inequality,& Math. Ann.& {\bf 246}
(1980) 249-250.

\ref Hi& M. Hilsum& Les espaces $L_p$ d'une alg\`{e}bre de von
Neumann d\'{e}finie par la d\'{e}riv\'{e}e spatiale,&J. Funct.
Anal.&{\bf 40} (1981) 151-169.

\ref I&H. Izumi&Construction of noncommutative $L_p$-spaces with a
complex parameter arising from modular actions,&Internat. J.
Math.&{\bf 8} (1997) 1029-1066.

\ref J& R. C. James& Uniformly non-square Banach Spaces,&Ann.
Math. &{\bf 80} (1964) 542--550.

\ref Ja&H. Jarchow& On weakly compact operators on
$C^*$-algebras,&Math. Ann. &{\bf 273} (1986) 341-343.

\ref JNRX& M. Junge, N.J. Nielsen, Z-J. Ruan, Q. Xu& ${\cal OL}_p$
spaces and the local structure of non-commutative $L_p$-spaces,&
Preprint&, 2001.

\ref JOR& M. Junge, N. Ozawa and J-Z. Ruan&On ${\cal OL}_{\infty}$
structure of nuclear $C^*$-algebras,& preprint.&

\ref Ju1&M. Junge& Embeddings of non-commutative $L_p$-spaces into
non-commutative $L\sb 1$-spaces, $1<p<2$,& Geom. Funct. Anal.&{\bf
10} (2000),  389-406.

\ref Ju2& M. Junge&  Doob's inequality for non-commutative
martingales,& J. Reine Angew. Math.& {\bf 549} (2002), 149-190.

\ref Ju3& M. Junge&  A Fubini type theorem for non-commutative
$L_p$ spaces,& preprint.&

\ref JX& M. Junge, Q. Xu& Non-commutative Burkholder/Rosenthal
inequalities and applications,& Annals of Proba.,& to appear.

\ref KaP& M. I. Kadec, A. Pe\l czy\'nski& Bases, lacunary
sequences and complemented subspaces in the spaces $L_p$, &Studia
Math.& {\bf 21} (1962), 161-176.

\ref KaR& R.V. Kadison, J.R. Ringrose&& Fundamentals of the theory
of operator algebras I and II,& Academic Press, 1983 and 1986.

\ref Ko1& H. Kosaki& Applications of the complex interpolation
method to a von Neumann Algebra: Non-commutative $L^p$-spaces,& J.
Funct. Anal.& {\bf 56} (1984), 29-78.

\ref Ko2&H. Kosaki&Applications of the uniform convexity of
noncommutative $L_p$-spaces,& Trans. Amer. Math. Soc.& {\bf 283}
(1984) 265-282.

\ref Ko3&H. Kosaki&On the continuity of the map
$\ph\to\vert\ph\vert$ from the predual of a von Neumann
Algebra,&J. Funct. Anal. & {\bf 59} (1984) 123-131.

\ref KrM &J. L. Krivine, B. Maurey&Espaces de Banach
stables,&Israel J. of Math.& {\bf 39} (1981) 273-295.

\ref LT& J. Lindenstrauss, L. Tzafriri&&Classical Banach spaces,
Vol. 1,&Springer Verlag, 1977.

\ref MS&V. D.  Milman, G. Schechtman&& Asymptotic theory of finite
dimensional normed spaces,& Lectures notes in Math. {\bf 1200},
Springer Verlag (1986).

\ref N& E. Nelson&Notes on non-commutative integration,&J. Funct.
Anal.&{\bf 15} (1974) 103-116.

\ref NO&P. W. Ng, N. Ozawa& A characterization of completely
1-complemented subspaces of non-commutative $L^1$-spaces,&
preprint.&

\ref O& N. Ozawa& Almost completely isometric embeddings between
preduals of von Neumann algebras, & J. Funct. Anal.& {\bf 186}
(2001), 329-341.

\ref P1& G. Pisier&& Non-commutative vector valued $L_p$-spaces
and completely $p$-summing maps, &Ast\'erisque {\bf 247}, (1998).

\ref P2& G. Pisier&& Introduction to  operator space
theory,&Cambridge Univ. Press, to appear.

\ref PX& G. Pisier, Q. Xu&Non-commutative martingale
inequalities,&
 Comm. Math. Phys.&{\bf 189} (1997) 667-698.

\ref R1& N. Randrianantoanina& Sequences in non-commutative
$L^p$-spaces,& Preprint,& 2000.

\ref R2& N. Randrianantoanina&  Embedding of $\ell_p$ in
non-commutative spaces,& Preprint,& 2000.

\ref R3& N. Randrianantoanina& Kadec-Pe\l czy\'nski decomposition
for Haagerup $L^p$-spaces,& Math. Camb. Philos. Soc.  &{\bf 132}
(2002), 137-154.

\ref R4& N. Randrianantoanina&  Non-commutative martingale
transforms,& J. Funct. Anal.& to appear.

\ref Ra& Y. Raynaud& On ultrapowers of non commutative
$L_{p}$-spaces,& J. Operator Theory,& to appear.

\ref RaX & Y. Raynaud, Q. Xu& On the structure of subspaces of
non-commutative $L_p$-spaces,&C. R. Acad. Sci. Paris&{\bf 333}
(2001) 213-218.

\ref Ro& H. P. Rosenthal& On subspaces of $L_p$,& Ann. Math.&{\bf
97} (1973) 344-373.

\ref S& F. A. Sukochev&  Non-isomorphism of $L_p$-spaces
associated with finite and infinite von Neumann algebras,&Proc.
Amer. Math. Soc.& {\bf 124} (1996), 1517-1527.

\ref SX& F. A. Sukochev, Q. Xu& Embedding of non-commutative
$L_p$-spaces, $p<1$,& Archiv der Math.& to appear.

\ref T& M. Takesaki&& Theory of operator algebras I,
&Springer-Verlag, New York-Heidelberg Berlin, 1979.

\ref Te1& M. Terp&& $L^p$-spaces associated with von Neumann
algebras,& Notes, Math. Institute, Copenhagen Univ. 1981.

\ref Te2& M. Terp& Interpolation spaces between a von Neumann
algebra and its predual,& J. Operator Th. &{\bf 29} (1977) 73-90.

\ref W& L. W. Weis& Banach lattices with the subsequence splitting
property,& Proc. Amer. Math. Soc. &{\bf 105} (1989) 87-96.

}

\vskip 0.5cm\parindent= 0pt {\eightsc Yves Raynaud: Equipe d'Analyse
Fonctionnelle, Institut de Math\'{e}matiques de Jussieu (CNRS), Case
186, 4, place Jussieu, 75252 Paris Cedex 05, France.}

{\eightsl E-mail:} {\eightpoint yr@ccr.jussieu.fr}

{\eightsc Quanhua Xu: Laboratoire de Math\'ematiques, Universit\'e de
Franche-Comt\'e, Route de Gray, 25030 Besan\c con cedex, France.}

{\eightsl E-mail:} {\eightpoint qx@math.univ-fcomte.fr}

\end